\documentclass[reqno,draft]{amsart}

\usepackage{amsmath}
\usepackage{amscd}
\usepackage{amssymb}
\usepackage{enumerate}
\usepackage{amsfonts}
\usepackage{graphicx}

\DeclareFontEncoding{OT2}{}{}

\numberwithin{equation}{section}

\newcommand{\rar}[1]{\stackrel{#1}{\longrightarrow}}

\newcommand{\be}{\beta}

\newcommand{\Ga}{\Gamma}
\newcommand{\de}{\delta}
\newcommand{\De}{\Delta}

\newcommand{\La}{\Lambda}

\newcommand{\sg}{\sigma}

\newcommand{\om}{\omega}
\newcommand{\Om}{\Omega}

\newcommand{\bC}{{\mathbb C}}

\newcommand{\bR}{{\mathbb R}}

\newcommand{\cB}{{\mathcal B}}

\newcommand{\cJ}{{\mathcal J}}
\newcommand{\cK}{{\mathcal K}}
\newcommand{\cL}{{\mathcal L}}

\newcommand{\cN}{{\mathcal N}}

\newcommand{\eq}[1]{(\ref{#1})}

\newcommand{\Ker}{\operatorname{Ker}}
\newcommand{\End}{\operatorname{End}}

\newcommand{\Aut}{\operatorname{Aut}}

\newcommand{\Wedge}{\bigwedge}

\renewcommand{\dot}{^{\bullet}}

\newcommand{\sbr}{\smallbreak}

\newcommand{\Ann}{\operatorname{Ann}}
\newcommand{\Image}{\operatorname{Im}}
\newcommand{\cliff}{\operatorname{Cliff}}
\newcommand{\Cou}{\mathcal{COU}}
\newcommand{\cou}[2]{\bigl[#1,#2\bigr]_{\mathrm{cou}}}

\newcommand{\Can}{\underline{\operatorname{Can}}}

\newcommand{\cJt}{\tilde{\cJ}}

\newcommand{\bra}{\langle}
\newcommand{\ket}{\rangle}
\newcommand{\pairing}{\bra\cdot,\cdot\ket}
\newcommand{\pair}[2]{\bra #1, #2\ket}
\newcommand{\eval}[2]{\bra #1\,\big\vert\,#2\ket}

\newcommand{\dimr}{\dim_{\bR}}
\newcommand{\dimc}{\dim_{\bC}}

\newcommand{\restr}[1]{\big\vert_{#1}}

\newcommand{\matr}[4]{\left(\begin{array}{cc} \! #1 & \! #2 \! \\ \! #3 & \! #4 \! \end{array}\right)}

\newtheorem{ithm}{Theorem}

\newtheorem{thm}{Theorem}[section]
\newtheorem{cor}[thm]{Corollary}
\newtheorem{lem}[thm]{Lemma}
\newtheorem{prop}[thm]{Proposition}

\theoremstyle{remark}
\newtheorem{rem}[thm]{Remark}

\newtheorem{example}[thm]{Example}
\newtheorem{examples}[thm]{Examples}
\newtheorem{defin}[thm]{Definition}

\title[Submanifolds of GCM]{Submanifolds of generalized (almost) complex manifolds}

\author[O.~Ben-Bassat and M.~Boyarchenko]{Oren Ben-Bassat\address{Oren Ben-Bassat:
Department of Mathematics, University of Pennsylvania,
Philadelphia, PA 19104-6395, e-mail: orenb@math.upenn.edu} \and
Mitya Boyarchenko\address{Mitya Boyarchenko: Department of
Mathematics, University of Chicago, Chicago, IL 60637, \\ e-mail:
mitya@math.uchicago.edu}}

\date{September 1, 2003}

\begin{document}

\begin{abstract}
The main goal of our paper is the study of several classes of
submanifolds of generalized complex manifolds. Along with the
generalized complex submanifolds defined by Gualtieri and Hitchin
in \cite{Gua}, \cite{H3} (we call these ``generalized Lagrangian
submanifolds'' in our paper), we introduce and study three other
classes of submanifolds. For generalized complex manifolds that
arise from complex (resp., symplectic) manifolds, all three
classes specialize to complex (resp., symplectic) submanifolds. In
general, however, all three classes are distinct. We discuss some
interesting features of our theory of submanifolds, and illustrate
them with a few nontrivial examples.

We then support our ``symplectic/Lagrangian viewpoint'' on the
submanifolds introduced in \cite{Gua}, \cite{H3} by defining the
``generalized complex category'', modelled on the constructions of
Guillemin-Sternberg \cite{GS} and Weinstein \cite{Wei2}. We argue
that our approach may be useful for the quantization of
generalized complex manifolds.
\end{abstract}

\maketitle

\tableofcontents

\section{Introduction}\label{s:intro}

\subsection{Motivation}\label{ss:motivation}

The goal of our paper is the study of certain special classes of
submanifolds of generalized complex manifolds. The notion of a
generalized complex manifold was introduced by Nigel Hitchin in
\cite{H2}. It contains as special cases both complex and
symplectic manifolds. Later Marco Gualtieri and Hitchin
(\cite{Gua}, \cite{H3}) have defined a class of submanifolds of
generalized complex manifolds which they have called ``generalized
complex submanifolds'' and which in this paper we call
``generalized Lagrangian submanifolds'', for the reasons explained
below. Hitchin's and Gualtieri's notion specializes to complex
submanifolds of complex manifolds and to Lagrangian submanifolds
of symplectic manifolds. In particular, if $M$ is a generalized
complex manifold, then, in general, neither $M$ itself nor the
points of $M$ are generalized complex submanifolds of $M$ in their
terminology.

\sbr

Our first task, then, was to find a different definition which
gives complex submanifolds in the complex case and symplectic
submanifolds in the symplectic case. It is also important to
understand the relationship between our notion of a generalized
complex submanifold and Hitchin and Gualtieri's notion; and, on
the other hand, to try to see if Hitchin and Gualtieri's notion
can be used in the theory of generalized complex manifolds in a
similar way that Lagrangian submanifolds are used in symplectic
geometry. For example, if $\phi:N\to M$ is a diffeomorphism
between two symplectic manifolds, it is well-known (see, e.g.,
\cite{dS}) that $\phi$ is a symplectomorphism if and only if the
graph of $\phi$ is a Lagrangian submanifold of $N\times M$ with
respect to the ``twisted product symplectic structure'' on
$N\times M$. It is natural to ask if the obvious generalization of
this result holds for generalized complex manifolds.

\sbr

More generally, suppose $(M,\om_M)$ is a symplectic manifold and
$N$ is a submanifold of $M$ which carries a symplectic form
$\om_N$. Again, it is easy to see that $\om_N=\om_M\big\vert_N$ if
and only if the graph of the inclusion map $N\hookrightarrow M$ is
an isotropic submanifold of $N\times M$ with respect to the
twisted product symplectic structure. On the other hand, if $M$ is
a real manifold, $N$ is a submanifold of $M$, and both $M$ and $N$
are equipped with complex structures, then $N$ is a complex
submanifold of $M$ if and only if the graph of the inclusion map
$N\hookrightarrow M$ is a complex submanifold of $N\times M$ with
respect to the product complex structure.

\sbr

A starting point for our discussion is the observation that the
conditions appearing in Gualtieri and Hitchin's definition of a
generalized complex submanifold can be relaxed naturally, so that
the resulting object specializes to complex submanifolds in the
complex case and to isotropic submanifolds in the symplectic case.
We call these objects ``generalized isotropic submanifolds.'' One
of the main goals of our work was trying to generalize the remarks
of the previous paragraph to find a relationship between
generalized isotropic submanifolds and our notion of generalized
complex submanifolds.

\sbr

Another goal was to define and study the ``generalized complex
category'', in the spirit of \cite{Wei}, \cite{Wei2}. The
definition is based on the notion of a generalized Lagrangian
submanifold. We believe that this ``category'' can be used to
study quantization of generalized complex manifolds. (The word
``category'' appears in quotes, because, as explained in
\cite{Wei} or \cite{Wei2}, a certain transversality condition has
to be satisfied for a pair of morphisms to be composable.)

\subsection{Main results}\label{ss:results} Most of our
paper is devoted to the study of generalized complex structures on
real vector spaces. Only at the end do we apply our results to
generalized complex manifolds. While we feel that our study of
subspaces of generalized complex vector spaces is more or less
complete, it is clear that our results for generalized complex
manifolds have a preliminary nature. Still, we believe that the
questions that we have raised are already fundamental at the
vector space level and thus have to be addressed. Furthermore, as
we show, the linear situation is far from being trivial (unlike in
the theory of complex or symplectic vector spaces), due to the
possibility of a $B$-field transformation.\footnote{These
transformations, along with the $\be$-field transformations, play
a special role in the theory. They were first introduced in
\cite{H2}. We recall their definition in \S\ref{ss:B-field}.} Some
unexpected features of the theory are described below.

\sbr

The first part of our paper is devoted to the study of induced
generalized complex structures on subspaces of generalized complex
vector spaces. The situation is quite similar to the one
considered by Theodore Courant in \cite{Cou}. In fact, many of our
basic constructions are straightforward modifications of the ones
introduced by Courant. However, there is an important difference:
whereas a Dirac structure on a real vector space induces a Dirac
structure on every subspace and quotient space, this is not so for
generalized complex structures. Namely, recall that if $V$ is a
finite-dimensional real vector space, then a generalized complex
structure on $V$ can be defined by a subspace $E\subset
V_{\bC}\oplus V_{\bC}^*$ which is isotropic with respect to the
quadratic form $Q(v\oplus f)=-f(v)$ and satisfies $V_{\bC}\oplus
V_{\bC}^*=E\oplus\overline{E}$, where the bar denotes the complex
conjugate. Now, given a real subspace $W\subseteq V$, we can
define, as in \cite{Cou}, the subspace $E_W\subset W_{\bC}\oplus
W_{\bC}^*$ consisting of all pairs of the form
$\bigl(w,f\big\vert_{W_{\bC}}\bigr)$, where $w\in W_{\bC}$ and
$f\in V_{\bC}$ are such that $(w,f)\in E$. It is clear that $E_W$
is still isotropic with respect to $Q$. However, the condition
$W_{\bC}\oplus W_{\bC}^*=E_W\oplus\overline{E}_W$ may no longer
hold.

\sbr

Thus we define $W$ to be a ``generalized complex subspace'' if the
last condition does hold, so that we get an induced generalized
complex structure. We then describe this induced structure in
terms of spinors (Proposition \ref{p:subspin}). The dual notion,
that of a quotient generalized complex structure, is also studied.
We introduce an operation which interchanges generalized complex
structures on a real vector space and on its dual; this operation
can then be used to reduce many statements about quotients to
statements about subspaces (and vice versa). Hence, naturally, we
devote more attention to subspaces than to quotients. We also
prove the ``classification theorem'' for generalized complex
vector spaces: every generalized complex structure is the
$B$-field transform of a direct sum of a complex structure and a
symplectic structure (Theorem \ref{t:classification}).

\sbr

Let us emphasize that there are several subtle points in the
theory of subspaces that we develop, which are not present in the
theory of complex or symplectic structures on vector spaces. For
example, if $W$ is a generalized complex subspace of a generalized
complex vector space $V$, then $V/W$ need {\em not} be a
generalized complex quotient of $V$ (cf. Example
\ref{e:subnotquot}). As another example of ``strange'' behavior,
we show the following. If $V$ is any generalized complex vector
space, there is a unique maximal subspace $S$ of $V$ such that $S$
is a GC subspace and the induced structure on $S$ is
$B$-symplectic (Theorem \ref{t:canon}). [Hereafter, we use the
expressions ``$B$-symplectic structure'' and ``$B$-complex
structure'' as shorthand for ``$B$-field transform of a symplectic
structure'' and ``$B$-field transform of a complex structure,''
respectively.] Furthermore, there is a canonical complex structure
$J$ on the quotient $V/S$. {\em However, $V/S$ may not be a
generalized complex quotient of $V$, in the sense of Definition
\ref{d:gcquot}.} (See Example \ref{ex:notquot}.) On the other
hand, if $V/S$ {\em is} a generalized complex quotient of $V$,
then the induced structure on $V/S$ is $B$-complex, and the
associated complex structure coincides with $J$.

\sbr

The interplay between $B$-field transforms and $\be$-field
transforms is also delicate. For instance, if a generalized
complex structure is both a $B$-field transform and a $\be$-field
transform of a (classical) complex structure, then it is, in fact,
complex. The corresponding statement for symplectic structures is
false. On the other hand, it is possible to start with a classical
symplectic structure, make a $B$-field transform, then a
$\be$-field transform, and arrive at a classical complex
structure. These features are discussed in some detail in
\S\ref{ss:symplB}.

\sbr

To describe the results of the second part of the paper, we need
to introduce more notation. Let $V$ be as above, and recall that
another way to describe a generalized complex structure on $V$ is
by giving an automorphism $\cJ_V$ of $V\oplus V^*$ which satisfies
$\cJ_V^2=-1$ and preserves the quadratic form $Q$ on $V\oplus
V^*$. [An equivalence with the previous definition is obtained by
assigning to such a $\cJ_V$ its $+i$-eigenspace $E$ in
$V_{\bC}\oplus V_{\bC}^*$.] For many considerations it is
convenient to write $\cJ_V$ in matrix form:
\[
\cJ_V=\matr{\cJ_1}{\cJ_2}{\cJ_3}{\cJ_4},
\]
where $\cJ_1:V\to V$, $\cJ_2:V^*\to V$, $\cJ_3:V\to V^*$ and
$\cJ_4:V^*\to V^*$ are linear maps. Let $W$ be a real subspace of
$V$. We say that $W$ is {\em split} if there exists a subspace
$N\subseteq V$ such that $V=W\oplus N$ and
$W\oplus\Ann(N)\subseteq V\oplus V^*$ is stable under $\cJ$. [In
fact, this is equivalent to the condition that there is a
decomposition $V=W\oplus N$ as a direct sum of generalized complex
structures, which explains the terminology.] On the other hand, we
define the ``twisted generalized complex structure'' on $V$ to be
the one corresponding to the automorphism
\[
\cJt_V=\matr{\cJ_1}{-\cJ_2}{-\cJ_3}{\cJ_4}.
\]
If $W$ is equipped with some generalized complex structure
$\cJ_W$, we say that this subspace {\em satisfies the graph
condition} (with respect to $\cJ_W$) if the graph of the inclusion
map $W\hookrightarrow V$ is a generalized isotropic subspace of
$W\oplus V$ with respect to the ``twisted product structure'' on
$W\oplus V$, which by definition is the direct sum of $\cJt_W$ and
$\cJ_V$.

\sbr

The following three results (see \S\ref{ss:subspaces}) explain the
relations between the three types of subspaces of a generalized
complex vector space $V$ that we have introduced.
\begin{ithm}
Every split subspace is a generalized complex subspace.
Furthermore, it satisfies the graph condition with respect to the
induced generalized complex structure.
\end{ithm}
\begin{ithm}
A generalized complex subspace $W\subseteq V$ satisfies the graph
condition with respect to the induced structure if and only if $W$
is invariant under $\cJ_1$.
\end{ithm}

\noindent It is worth emphasizing that if $W$ is a generalized
complex subspace and it satisfies the graph condition with respect
to {\em some} generalized complex structure, this structure does
not necessarily have to be the induced one. More precisely:
\begin{ithm}
Suppose that $W\subseteq V$ is an arbitrary subspace and $\cJ_W$
is a generalized complex structure on $W$ such that the graph
condition is satisfied. If $W$ is a generalized complex subspace,
then the graph condition is also satisfied with respect to the
induced generalized complex structure, and $\cJ_W$ is a
$\be$-field transform of the induced structure on $W$.
\end{ithm}
In Remark \ref{r:graphnotsub}, we show that a subspace satisfying
the graph condition with respect to {\em some} generalized complex
structure does not have to be a generalized complex subspace in
our sense. This may be seen as evidence of the fact that our
definition of the graph condition is too weak. One possible
explanation for this is that, unlike the classical symplectic
situation, a generalized isotropic subspace which has half the
dimension of the ambient space is not necessarily generalized
Lagrangian. Thus, for example, if $\mu:V\to W$ is an isomorphism
of real vector spaces, and $V$, $W$ are equipped with generalized
complex structures such that the graph of $\mu$ is generalized
isotropic with respect to the twisted product structure, then it
does not necessarily follow that $\mu$ induces an isomorphism
between the two structures. On the other hand, we have the
following
\begin{ithm}\label{it:4}
Let $V$, $W$ be generalized complex vector spaces, and let
$\mu:V\to W$ be an isomorphism of the underlying real vector
spaces. If the graph of $\mu$ is generalized Lagrangian with
respect to the twisted product structure on $V\oplus W$, then
$\mu$ induces an isomorphism between the two generalized complex
structures.
\end{ithm}

\noindent The global version of this result is given in Theorem
\ref{t:graphiso}.

\sbr

In the second part of the paper we generalize all the notions
introduced previously to generalized complex manifolds. In
addition, we study the analogue of the notion of an ``admissible
function'' introduced by Courant \cite{Cou}. This leads to the
definition of the {\em Courant sheaf}: if $M$ is a generalized
complex manifold, with $\cJ_M$ the corresponding automorphism of
$TM\oplus T^*M$ and $\overline{E}$ the $-i$ eigenbundle of $\cJ_M$
on $T_{\bC}M\oplus T_{\bC}^*M$, the sections of the Courant sheaf
are pairs $(f,X)$, where $f$ is a complex-valued $C^\infty$
function and $X$ is a section of $T_{\bC}M$ such that
$df+X\in\overline{E}$. We prove that the Courant sheaf is a sheaf
of local rings and also a sheaf of Poisson algebras. Moreover, for
a split submanifold $N$ of a generalized complex manifold $M$, we
define a natural operation of ``pullback'' of sections of the
Courant sheaf on $M$ to those on $N$ (Theorem \ref{t:pullcour}).
The generalized complex structure on $M$ can be recovered from the
Courant sheaf together with its embedding into $C^\infty_M\oplus
T_{\bC}M$ in case $M$ is complex or symplectic, but it seems that,
in general, additional information is required.

\sbr

The last part of the paper is devoted to the generalized complex
``category''. Recall that Weinstein \cite{Wei2} has introduced the
symplectic ``category'', based on some ideas of Guillemin and
Sternberg \cite{GS}. The objects of this ``category'' are
symplectic manifolds. If $M$ and $N$ are two symplectic manifolds,
the set of morphisms from $M$ to $N$ is defined to be the set of
submanifolds of $M\times N$ that are Lagrangian with respect to
the ``twisted product symplectic structure'' on $M\times N$. Two
morphisms, $M\to N$ and $N\to K$, can be composed provided a
certain transversality condition is satisfied (see \cite{Wei2} for
details). In \S\ref{ss:category} we generalize Weinstein's
construction by replacing symplectic manifolds with generalized
complex manifolds, and Lagrangian submanifolds with generalized
Lagrangian submanifolds. In the linear case (\S\ref{ss:catgcvs}),
where the transversality condition disappears, we get an honest
category, thereby generalizing a construction of Guillemin and
Sternberg \cite{GS}. We would like to point out that one of the
applications of Guillemin-Sternberg-Weinstein construction is to
deformation quantization and geometric quantization of symplectic
manifolds (see \cite{GS}, where the quantization problem is solved
completely in the linear case, and \cite{Wei}, \cite{Wei2} for a
discussion of the global case). We hope that, similarly, the
constructions of \S\ref{s:category} can be used for quantization
of generalized complex structures. We plan to study this in detail
in subsequent publications.

\subsection{Structure of the paper}\label{ss:structure} Our paper
is organized as follows. In \S\ref{s:linalg} we recall the basic
definitions of the theory of generalized complex vector spaces. We
follow \cite{Gua} and \cite{H3} for the most part. We also list
the basic definitions and results from the theory of spinors,
referring to the book \cite{Che} for proofs. In \S\ref{s:subquot}
we introduce and study the notions of ``generalized complex
subspaces'' and ``generalized complex quotients,'' and prove the
classification theorem, as well as certain weaker but more
``canonical'' analogues of the theorem. In \S\ref{s:split} we
introduce the notions of ``split subspaces'' and ``subspaces
satisfying the graph condition'' and describe in detail the
relationship of these notions with the more general notion of a
generalized complex subspace.

\sbr

In \S\ref{s:gcm} we recall the definitions of generalized almost
complex manifolds and generalized complex manifolds. For the
latter, we study the natural generalization of the notion of
``admissible function'' defined by Courant, and we introduce the
``Courant sheaf.'' We then define the global analogues of
generalized complex subspaces and split subspaces: ``generalized
complex submanifolds'' and ``split submanifolds'' of generalized
complex manifolds. We end the section by defining the ``pullback''
of sections of the Courant sheaf on a generalized complex manifold
to a split submanifold. In \S\ref{s:category}, we study the
``graph condition'' for submanifolds of generalized complex
manifolds. This idea allows us to define the ``generalized complex
category'' \`a la Weinstein \cite{Wei2}. We finish by studying the
simpler category of generalized complex vector spaces, modelled on
the construction of Guillemin and Sternberg \cite{GS}.

\sbr

In \S\ref{s:examples} we illustrate our main constructions and
results with two types of examples: $B$-field transforms of
complex and symplectic structures. (A few scattered examples also
appear throughout the text.) More specifically, we characterize
the GC subspaces (resp. subspaces satisfying the graph condition,
resp. split subspaces) of $B$-complex and $B$-symplectic vector
spaces, we characterize generalized isotropic (resp. coisotropic,
resp. Lagrangian) subspaces of complex and symplectic vector
spaces, and we state Gualtieri's characterization of generalized
Lagrangian submanifolds of $B$-complex and $B$-symplectic
manifolds. We then explain some interesting features of $B$-field
transforms of symplectic manifolds. Finally, we discuss in detail
some nontrivial counterexamples that we mentioned in
\S\ref{ss:results}.

\sbr

We have only included the proofs of the most trivial results in
the main body of the paper. All other proofs appear in a separate
section (\S\ref{s:proofs}) at the end of the article. This
clarifies the exposition, and allows us to present our results in
not necessarily the same order in which they are proved.

\sbr

For the reader's convenience, a more detailed summary appears at
the beginning of each section.

\subsection{Acknowledgements}\label{ss:thanks} We are grateful
to Marco Gualtieri for sharing some of his notes with us, prior to
the publication of his thesis. In particular, several definitions
and constructions of \S\S\ref{s:linalg}, \ref{s:gcm},
\ref{s:examples} are based on these notes. We would also like to
thank Tony Pantev for the suggestion to work on generalized
complex manifolds, as well as for many useful comments on our
paper and suggestions for improvement.

\subsection{Notation and conventions}\label{ss:notation}
All vector spaces considered in this paper will be
finite-dimensional. We will be dealing with both real and complex
vector spaces. If $V$ is a real vector space, we denote by
$V_{\bC}$ its complexification. If $W\subseteq V$ is a real
subspace, we will identify $W_{\bC}$ with a subspace of $V_{\bC}$
in the usual way. The complex conjugation on $\bC$ extends to an
$\bR$-linear automorphism of $V_{\bC}$ which we denote by
$x\mapsto \bar{x}$. Recall that a complex subspace $U$ of
$V_{\bC}$ is of the form $W_{\bC}$ for some real subspace
$W\subseteq V$ if and only if $U=\bar{U}$.

\sbr

The dual of a real or complex vector space $V$ will be denoted by
$V^*$. We will always implicitly identify $V$ with $(V^*)^*$. If
$W\subseteq V$ is a real (resp., complex) subspace, we write
$\Ann(W)=\bigl\{f\in V^*\Big\lvert f\big\vert_W\equiv 0\bigr\}$.
Similarly, if $U\subseteq V^*$ is a real or complex subspace, we
write $\Ann(U)=\bigl\{ v\in V \big\vert f(v)=0\ \forall\,f\in U
\bigr\}$.

\sbr

If $V$ is real, we identify $(V_{\bC})^*$ with $(V^*)_{\bC}$ in
the natural way. The natural projections $V\oplus V^*\to V$ and
$V\oplus V^*\to V^*$, or their complexifications, will be denoted
by $\rho$ and $\rho^*$, respectively. We have the standard
indefinite nondegenerate symmetric bilinear pairing $\pairing$ on
$V\oplus V^*$, given by\footnote{Our convention here follows that
of \cite{H3}.}
\[
\pair{v+f}{w+g}=-\frac{1}{2}\bigl(f(w)+g(v)\bigr) \text{  for all
} v,w\in V,\  f,g\in V^*.
\]
We use the same notation for the complexification of $\pairing$.
Usually, if $v\in V$ and $f\in V^*$, we will write either
$\eval{f}{v}$ or $\eval{v}{f}$ for $f(v)$. The pairing $\pairing$
corresponds to the quadratic form $Q(v+f)=-f(v)$.

\sbr

The exterior algebra of a vector space $V$ will be denoted by
$\Wedge\dot V$. If $x\in V$ and $\xi\in\Wedge\dot V^*$, we write
$\iota_x\xi$ for the contraction of $x$ with $\xi$. An element
$B\in\Wedge^2 V^*$ will be identified with the linear map $V\to
V^*$ that it induces; similarly, an element $\be\in\Wedge^2 V$
will be identified with the induced linear map $V^*\to V$. We will
frequently consider endomorphisms of $V\oplus V^*$. These will
usually be written as $2\times 2$ matrices
\[
T=\matr{T_1}{T_2}{T_3}{T_4},
\]
with the understanding that $T_1:V\to V$, $T_2:V^*\to V$,
$T_3:V\to V^*$ and $T_4:V^*\to V^*$ are linear maps. In this
situation, elements of $V\oplus V^*$ will be written as column
vectors: $(v,f)^t$, where $v\in V$, $f\in V^*$.

\sbr

For us, a {\em real} (resp., {\em complex}) {\em manifold} will
always mean a $C^\infty$ real (resp., complex-analytic)
finite-dimensional manifold. If $M$ is a real manifold, by
$C^\infty_M$ we denote the sheaf of (germs of) {\em
complex-valued} $C^\infty$ functions on $M$. We will denote by
$TM$ and $T^*M$ the tangent and cotangent bundles of $M$, and by
$T_{\bC}M$ and $T^*_{\bC}M$ their complexifications. The pairing
$\pairing$ defined above at the level of vector spaces extends
immediately to a pairing $\pairing:V_{\bC}\times V^*_{\bC}\to
C^\infty_M$, where $V$ is any vector bundle on $M$ and $V_{\bC}$
is its complexification. The de Rham differential on $\Wedge\dot
T^*_{\bC}M$ will be denoted by $\om\mapsto d\om$. For a section
$X$ of $T_{\bC}M$, we will denote by $\iota_X$ and $\cL_X$ the
operators of contraction with $X$ and the Lie derivative in the
direction of $X$, which are $\bC$-linear operators on $\Wedge\dot
T^*_{\bC}M$.

\sbr

By a {\em complex structure} on a real vector space $V$ we mean an
$\bR$-linear automorphism $J$ of $V$ satisfying $J^2=-1$. (Other
authors call such a $J$ an {\em almost complex structure}.) The
following abbreviations will be frequently used in our paper:
GC=``generalized complex'', GAC=``generalized almost complex'',
GCS=``GC structure'', GACS=``GAC structure'', GCM=``generalized
complex manifold'', GACM=``generalized almost complex manifold'',
GCY=``generalized Calabi-Yau manifold''.

\section{Generalized complex linear algebra}\label{s:linalg}

\subsection{Summary} This section describes the
various ways of defining a ``generalized complex structure'' on a
real vector space $V$. Our main purpose here is to fix the
notation used in subsequent arguments, and to make our paper as
self-contained as possible. The definition of a GCS (Generalized
Complex Structure) in terms of the subspace $E\subseteq
V_{\bC}\oplus V_{\bC}^*$ first appeared in \cite{H2}, the
definition in terms of an $\bR$-linear automorphism $\cJ$ of
$V\oplus V^*$ was introduced in \cite{Gua}, and the definition in
terms of spinors appears implicitly in \cite{H2} and explicitly in
\cite{H3}. All these definitions are studied in detail in
\cite{Gua}.

\sbr

Next we define an operation which interchanges generalized complex
structures on a real vector space with those on its dual. While
being completely straightforward, this operation is rather useful
in the study of sub- and quotient GC structures (see
\S\ref{s:subquot}).

\sbr

Finally, following \cite{Gua} and \cite{H3}, we introduce
$B$-field and $\be$-field transforms of GC structures. We describe
$B$-field transforms in terms of the three main definitions of GC
structures. We also recall how complex and symplectic structures
give rise to GC structures (\cite{Gua}), and we give a simple
characterization of the GC structures which are $B$-field (resp.,
$\be$-field) transforms of complex or symplectic structures.

\subsection{Equivalent definitions}\label{ss:defgcvs} We begin by
giving three ways of defining a GCS on a real vector space, and
then showing their equivalence. We will see later that in various
contexts, one of these descriptions may be better than the others,
but there is no description that is ``the best'' universally. Thus
we prefer to work with all of them interchangeably.
\begin{prop}\label{p:eqvdef}
Let $V$ be a real vector space. There are natural bijections
between the following $3$ types of structure on $V$.
\begin{enumerate}[(1)]
\item A subspace $E\subset V_{\bC}\oplus V_{\bC}^*$ such that
$E\cap\overline{E}=(0)$ and $E$ is maximally isotropic with
respect to the standard pairing $\pairing$ on $V_{\bC}\oplus
V_{\bC}^*$.
\item An $\bR$-linear automorphism $\cJ$ of $V\oplus V^*$ such
that $\cJ^2=-1$ and $\cJ$ is orthogonal with respect to
$\pairing$.
\item A pure spinor $\phi\in\Wedge\dot V^*_{\bC}$, defined up to
multiplication by a nonzero scalar, such that
$\pair{\phi}{\bar{\phi}}_M\neq 0$, where $\pairing_M$ denotes the
Mukai pairing $($\cite{Che}, \cite{H2}, \cite{H3}$)$.
\end{enumerate}
\end{prop}
The terminology (3) is explained in \S\ref{ss:spinors} below, and
the bijections are made explicit in \S\ref{ss:proofslinalg}. A
structure of one of these types will be called a {\em generalized
complex structure} on $V$, and a real vector space equipped with a
GC structure will be called a {\em generalized complex vector
space}.
\begin{rem}\label{r:evendim}
In \cite{H1}, \cite{H2} it is assumed from the very beginning that
generalized complex manifolds have even dimension (as real
manifolds). The reader may also assume that our GC vector spaces
have even real dimension; this is used in an essential way in the
arguments that involve spinors (cf. Theorem
\ref{t:spinors}(\ref{i:pair})). However, we can prove that, with
the definitions of a GC structure on $V$ in terms of $E\subset
V_{\bC}\oplus V_{\bC}^*$ or $\cJ\in\Aut_{\bR}(V\oplus V^*)$, the
even-dimensionality of $V$ is automatic (Corollary
\ref{c:evendim}). We invite the reader to check that our proof of
Corollary \ref{c:evendim} does not depend on any statements about
spinors. The result was also proved in \cite{H3} by a different
method (by showing that a GC structure on a real manifold induces
an almost complex structure).
\end{rem}

\subsection{Recollection of spinors}\label{ss:spinors}
Here we list the basic definitions and results from the theory of
spinors, following \cite{Che} (and, in a few places, \cite{H3}),
that will be used in the sequel. Chevalley's description of the
theory is somewhat more general than presented below; however, the
reader should have no difficulty modifying the notation and
results of \cite{Che} to fit our situation.

\sbr

Let $V$ be a real vector space, and write $V^d=V\oplus V^*$.
Recall the quadratic form $Q(v+f)=-\eval{v}{f}=-\iota_v(f)$ on
$V^d$. We let $\cliff(V)$ denote the {\em Clifford algebra} of the
complexification $V^d_{\bC}$ corresponding to the quadratic form
$Q$. This means (see \cite{Che}) that $\cliff(V)$ is the quotient
of the tensor algebra on $V^d_{\bC}$ by the ideal generated by all
expressions of the form $e-Q(e)$, $e\in V^d_{\bC}$. Note that if
$E$ is any isotropic (w.r.t. $Q$) subspace of $V^d_{\bC}$, then
the subalgebra of $\cliff(V)$ generated by $E$ is naturally
isomorphic to the exterior algebra $\Wedge\dot E$. In particular,
$\cliff(V)$ contains $S(V):=\Wedge\dot V^*_{\bC}$ in a natural
way. We call $S(V)$ the space of {\em spinors}. The subspaces of
$S(V)$ spanned by forms of even and odd degree are denoted by
$S_+(V)$ and $S_-(V)$, respectively. We call $S_+(V)$ (resp.,
$S_-(V)$) the space of {\em even} (resp., {\em odd}) {\em
half-spinors}.

\sbr

There exists a representation of $\cliff(V)$ on $S(V)$, with the
action of the generating subspace $V^d_{\bC}\subset\cliff(V)$
being defined by
\[
(v+f)\cdot\phi=\iota_v\phi+f\wedge\phi,\ \ v\in V_{\bC},\,f\in
V^*_{\bC},\,\phi\in S(V).
\]
Let $\phi\in S(V)$, $\phi\neq 0$, and let $E$ be the annihilator
of $\phi$ in $V^d_{\bC}$, with respect to the action defined
above. Observe that if $e\in E$, then $0=e\cdot
e\cdot\phi=Q(e)\cdot\phi$, whence $Q(e)=0$. Thus
$Q\big\vert_E\equiv 0$, so $E$ is isotropic with respect to $Q$.
\begin{defin}\label{d:puresp}
We say that $\phi$ is a {\em pure spinor} if $E$ is maximally
isotropic, i.e., $\dim_{\bC} E=\dim_{\bR} V$. In this case, we
also call $\phi$ a {\em representative spinor} for $E$.
\end{defin}
Recall from \S\ref{ss:notation} that complex conjugation acts on
$V^d_{\bC}$ and on $S(V)$. It is clear that $\phi\in S(V)$ is a
pure spinor if and only if $\bar{\phi}$ is such; and if $\phi$ is
a representative spinor for a maximally isotropic subspace
$E\subset V^d_{\bC}$, then $\bar{\phi}$ is a representative spinor
for $\overline{E}$, and vice versa.

\sbr

The last ingredient that we need is the {\em Mukai pairing}
$\pairing_M:S(V)\times S(V)\to \Wedge^n V$, where $n=\dim_{\bR} V$
(see \cite{Che}, \S{}III.3.2; \cite{H2}, \cite{H3}). However, we
have decided to omit the definition of this pairing, since it is
described in detail in \cite{Che} and \cite{H3}. The only fact
that we need about this pairing is Theorem
\ref{t:spinors}(\ref{i:pair}), which is proved in \cite{H3}. We
refer the reader to \cite{H3} for a detailed discussion of the
Mukai pairing that is most relevant to our situation.

\sbr

We can now state the main results from the theory of spinors that
we will use in the paper.
\begin{thm}\label{t:spinors}
\begin{enumerate}[(a)]
\item\label{i:evodd} Every pure spinor is either even or odd: if $\phi\in S(V)$
is pure, then $\phi\in S_+(V)$ or $\phi\in S_-(V)$.
\item\label{i:ann} Let $E\subset V^d_{\bC}$ be a maximally isotropic subspace;
then the space
\[
A_E=\bigl\{\phi\in S(V)\big\vert e\cdot\phi=0\ \forall\,e\in
E\bigr\}
\]
is one dimensional. Thus the nonzero elements of $A_E$ are
precisely the representative spinors for $E$. In other words, we
have a bijective correspondence between maximally isotropic
subspaces of $V^d_{\bC}$ and pure spinors modulo nonzero scalars.
\item\label{i:inters} Let $E,F\subset V^d_{\bC}$ be maximally isotropic subspaces
corresponding to pure spinors $\phi,\psi\in S_{\pm}(V)$. Then
$E\cap F=(0)$ if and only if $\pair{\phi}{\psi}_M\neq 0$.
\item\label{i:conj} In particular, if $E\subset V^d_{\bC}$ is a maximally isotropic
subspace corresponding to a pure spinor $\phi$, then
$E\cap\overline{E}=(0)$ if and only if
$\pair{\phi}{\bar{\phi}}_M\neq 0$.
\item\label{i:pure} A spinor $\phi\in S(V)$ is pure if and only if it can be
written as
\begin{equation}\label{e:puresp}
\phi=c\cdot\exp(u)\wedge f_1\wedge\dotsb\wedge f_k,
\end{equation}
where $c\in\bC^\times$, $u\in\Wedge^2 V^*$, and $f_1,\dotsc,f_k$
are linearly independent elements of $V^*_{\bC}$. If $\phi$ is
representative for a maximally isotropic subspace $E\subset
V^d_{\bC}$, then $f_1,\dotsc,f_k$ form a basis for $E\cap
V^*_{\bC}$.
\item\label{i:pair} Suppose that $\dim_{\bR} V$ is even,
and let $\phi\in S(V)$ be a pure spinor written in the form
\eq{e:puresp}. Then, up to a nonzero scalar multiple, we have
\[
\pair{\phi}{\bar{\phi}}_M=(u-\bar{u})^p\wedge
f_1\wedge\dotsb\wedge
f_k\wedge\bar{f}_1\wedge\dotsb\wedge\bar{f}_k,
\]
where $p=n/2-k$.
\end{enumerate}
\end{thm}
\begin{proof}[Proofs]
(a) \cite{Che}, III.1.5. (b) The proof is contained in \cite{Che},
\S{}III.1. (c) \cite{Che}, III.2.4. (d) This is immediate from
(c). (e) \cite{Che}, III.1.9. (f) See \cite{H3}.
\end{proof}

\subsection{Duality}\label{ss:duality} Given a real vector space $V$,
let $\tau:V\oplus V^*\to V^*\oplus V$ be the transposition of the
two summands. Recalling the natural identification of $V^{**}$
with $V$, we can also view $\tau$ as an isomorphism between
$V\oplus V^*$ and $V^*\oplus V^{**}$. We will continue to denote
by $\tau$ the induced isomorphism $V_{\bC}\oplus V^*_{\bC}\to
V^*_{\bC}\oplus V^{**}_{\bC}\cong V^*_{\bC}\oplus V_{\bC}$.
\begin{prop}
The isomorphism $\tau$ induces a bijection between GC structures
on $V$ and GC structures on $V^*$. Explicitly, $E\subset
V_{\bC}\oplus V^*_{\bC}$ is a GCS on $V$ if and only if $\tau(E)$
is a GCS on $V^*$. If $E$ corresponds to $\cJ\in\Aut_{\bR}(V\oplus
V^*)$ as in Proposition \ref{p:eqvdef}, then $\tau(E)$ corresponds
to $\tau\circ\cJ\circ\tau^{-1}$.
\end{prop}
\begin{proof}
If is clear that $\tau$ preserves the form $\pairing$, and if
$E\subseteq V_{\bC}\oplus V^*_{\bC}$ is any subspace, then
$\tau(\overline{E})=\overline{\tau(E)}$. The proposition follows
trivially from these observations.
\end{proof}

\subsection{$B$- and $\be$-field transforms}\label{ss:B-field}
Let $V$ be a real vector space and $B\in\Wedge^2 V^*$. We form the
matrix (see \S\ref{ss:notation})
\[
\cB:=\matr{1}{0}{B}{1}.
\]
It is straightforward to verify that $\cB$ is an orthogonal
automorphism of $V\oplus V^*$. Thus, if $E$ is a GCS on $V$, then
$\cB\cdot E$ is another one. We call $\cB\cdot E$ the {\em
$B$-field transform of $E$ defined by $B$}. Similarly, if
$\be\in\Wedge^2 V$, then the matrix
\[
\cB':=\matr{1}{\be}{0}{1}
\]
also acts on GC structures on $V$; if $E$ is one, then $\cB'\cdot
E$ will be called the {\em $\be$-field transform of $E$ defined by
$\be$}. If we look at GC structures in terms of the corresponding
orthogonal automorphisms $\cJ$ of $V\oplus V^*$, then the actions
of $B$ and $\be$ are given by $\cJ\mapsto\cB\cdot\cJ\cdot\cB^{-1}$
and $\cJ\mapsto\cB'\cdot\cJ\cdot\cB'^{-1}$, respectively. We can
also describe $B$-field transforms in terms of spinors:
\begin{prop}\label{p:Bspin}
Suppose that a GCS on a real vector space $V$ is defined by a pure
spinor $\phi\in\Wedge\dot V_{\bC}^*$, and let $B\in\Wedge^2 V^*$.
Then the $B$-field transform of this structure corresponds to the
pure spinor $\exp(-B)\wedge\phi$.
\end{prop}
We leave to the reader the task of describing $\be$-field
transforms in a similar way.
\begin{rem}\label{r:dualityBbeta}
Suppose that $E$ is a GCS on a real vector space $V$ and $E'$ is
the $B$-field transform of $E$ defined by $B\in\Wedge^2 V^*$.
Then, obviously, $\tau(E')$ is the $\be$-field transform of
$\tau(E)$, defined by the same $B\in\Wedge^2 V^*$ (but viewed now
as a bivector on $V^*$). Thus, the operation $\tau$ interchanges
$B$- and $\be$-field transforms.
\end{rem}
\begin{defin}\label{d:complsympl}
Let $V$ be a real vector space.
\begin{enumerate}[(1)]
\item Let $J$ be a (usual) complex structure on $V$. Then
\[
\cJ=\matr{J}{0}{0}{-J^*}
\]
is a GC structure on $V$. If $\cJ$ is a GCS on $V$ that can be
written in this form, we say that $\cJ$ is {\em complex}.
\item A $B$-field (resp., $\be$-field) transform of a complex GCS on $V$
will be referred to as a {\em $B$-complex} (resp., {\em
$\be$-complex}) structure on $V$.
\item Let $\om$ be a symplectic form on $V$ (i.e.,
a nondegenerate form $\om\in\Wedge^2 V^*$). Then
\[
\cJ=\matr{0}{-\om^{-1}}{\om}{0}
\]
is a GC structure on $V$. If $\cJ$ is a GCS on $V$ that can be
written in this form, we say that $\cJ$ is {\em symplectic}.
\item A $B$-field (resp., $\be$-field) transform of a symplectic
GCS on $V$ will be referred to as a {\em $B$-symplectic} (resp.,
{\em $\be$-symplectic}) structure on $V$.
\end{enumerate}
\end{defin}

\begin{prop}\label{p:characterization}
Let $V$ be a real vector space and $E\subset V_{\bC}\oplus
V_{\bC}^*$ a GCS on $V$. Write
\[
\cJ=\matr{\cJ_1}{\cJ_2}{\cJ_3}{\cJ_4}
\]
for the corresponding orthogonal automorphism of $V\oplus V^*$.
Then
\begin{enumerate}[(a)]
\item $E$ is $B$-complex $($resp., $\be$-complex$)$
$\iff$ $\rho(E)\cap\rho(\overline{E})=(0)$ $\bigl($resp.,
$\rho^*(E)\cap\rho^*(\overline{E})=(0)${}$\bigr)$ $\iff$ $\cJ_2=0$
$($resp., $\cJ_3=0${}$)$ $\iff$ $V^*_{\bC}=(V^*_{\bC}\cap
E)+(V^*_{\bC}\cap\overline{E})$ $\bigl($resp.,
$V_{\bC}=(V_{\bC}\cap E)+(V_{\bC}\cap\overline{E})${}$\bigr)$;
\item $E$ is complex if and only if it is both $B$-complex and
$\be$-complex;
\item $E$ is $B$-symplectic $($resp., $\be$-symplectic$)$ $\iff$
$E\cap V_{\bC}^*=(0)$ $\bigl($resp., $E\cap V_{\bC}=(0)${}$\bigr)$
$\iff$ $\cJ_2$ is an isomorphism $($resp., $\cJ_3$ is an
isomorphism$)$ $\iff$ $\rho(E)=V_{\bC}$ $($resp.,
$\rho^*(E)=V_{\bC}^*${}$)$;
\item $E$ is symplectic if and only if $\cJ_1=0$.
\end{enumerate}
\end{prop}
The analogue of part (b) fails in the symplectic case, as we will
see in \S\ref{ss:symplB}.

\section{Generalized complex subspaces and quotients}\label{s:subquot}

\subsection{Summary} In this section, we modify the construction of
\cite{Cou} to define the notion of a generalized complex subspace
of a GC vector space. The main difference with the situation in
\cite{Cou} is that, while a Dirac structure on a real vector space
induces a Dirac structure on each of its subspaces, this is not so
for generalized complex structures. We describe the induced GC
structure in terms of spinors, and study the behavior of GC
subspaces under a $B$-field transform of the ambient structure.

\sbr

We then define and study the dual notion of a quotient GC
structure. Using the ``duality operation'' $\tau$ introduced in
\S\ref{ss:duality}, we show that the study of GC quotients is
essentially equivalent to the study of GC subspaces. We also
mention a nontrivial unexpected counterexample.

\sbr

The last subsection is devoted to the problem of
``classification'' of generalized complex structures. It turns out
that the $B$-complex and $B$-symplectic structures are fundamental
examples of GC structures, in the sense that an arbitrary GC
structure is a $B$-field transform of a direct sum of a complex
and a symplectic structure (Theorem \ref{t:classification}). A
consequence of our classification theory is the fact that every GC
vector space has even dimension over $\bR$.

\subsection{Definition of GC subspaces}\label{ss:gcsubspace}
Let $V$ be a real vector space with a GC structure $E\subset
V_{\bC}\oplus V^*_{\bC}$. Given a subspace $W\subseteq V$, we
define (following \cite{Cou})
\begin{equation}\label{e:ew}
E_W=\Bigl\{ \bigl(\rho(e),\rho^*(e)\big\vert_{W_{\bC}}\bigr)
\Big\vert e\in E\cap \bigl(W_{\bC}\oplus V_{\bC}^*\bigr) \Bigr\}.
\end{equation}
Clearly, $E_W$ is an isotropic subspace of $W_{\bC}\oplus
W_{\bC}^*$. It is straightforward to compute that $\dim_{\bC}
E_W=\dim_{\bR} W$ (see Lemma \ref{l:dim}); thus $E_W$ is, in fact,
maximally isotropic.
\begin{defin}\label{d:gcsub}
We say that $W$ is a {\em generalized complex subspace} of $V$ if
$E_W\cap \overline{E}_W=(0)$, so that $E_W$ is a generalized
complex structure on $W$. In this case, we denote by $\cJ_W$ the
orthogonal automorphism of $W\oplus W^*$ corresponding to $E_W$.
\end{defin}
We note that a Dirac structure on a real vector space $V$ induces
a Dirac structure on every subspace of $V$ (see \cite{Cou}). In
our situation, however, not every subspace of a GC vector space is
a GC subspace. The reason is the extra condition
$E_W\cap\overline{E}_W=(0)$, which has no counterpart in the
theory of Dirac structures.
\begin{examples}
\begin{enumerate}[(1)]
\item Consider a complex GCS on $V$, corresponding to a complex structure
$J$ on $V$. Then a real subspace $W\subseteq V$ is a GC subspace if and
only if $W$ is stable under $J$. In this case, the induced GC structure
on $W$ is also complex, corresponding to $J\big\vert_W$.
\item Consider a symplectic GCS on $V$, corresponding to a symplectic
form $\om$ on $V$. Then a real subspace $W\subseteq V$ is a GC
subspace if and only if $\om\big\vert_W$ is nondegenerate on $W$.
In this case, the induced GC structure on $W$ is also symplectic,
corresponding to $\om\big\vert_W$.
\end{enumerate}
\end{examples}
This follows from Proposition \ref{p:subBcs}, where we give a
complete description of GC subspaces of $B$-complex and
$B$-symplectic vector spaces.
\begin{rem}
In general, it seems that if $W$ is a GC subspace of a GC vector
space $V$, then there is no simple relationship between $\cJ_W$
and $\cJ_V$. Later we will see that such a relationship exists for
special types of GC subspaces.
\end{rem}
It is important to understand the passage from a GCS on $V$ to the
induced GCS on a GC subspace $W\subseteq V$ in terms of spinors.
To that end, we have the following
\begin{prop}\label{p:subspin}
Let $V$ be a real vector space with a GC structure $E\subset
V_{\bC}\oplus V_{\bC}^*$. Let $W\subseteq V$ be any subspace, and
let $j:W\hookrightarrow V$ denote the inclusion map. Then there
exists a representative spinor for $E$ of the form
$\phi=\exp(u)\wedge f_1\wedge\dotsb\wedge f_k$, such that, for
some $1\leq l\leq k$, $j^*(f_1),\dotsc,j^*(f_l)$ are a basis of
$\Ann\bigl(\rho(E)\cap W_{\bC}\bigr)\subseteq W_{\bC}^*$, and
$f_{l+1},\dotsc,f_k$ are a basis of
$\Ann\bigl(\rho(E)+W_{\bC}\bigr)$. Moreover, $\phi_W:=\exp(j^*
u)\wedge j^*(f_1)\wedge\dotsb\wedge j^*(f_l)$ is a representative
spinor for $E_W$. In particular, $W$ is a GC subspace of $V$ if
and only if $\pair{\phi_W}{\bar{\phi}_W}_M\neq 0$.
\end{prop}
\begin{cor}\label{c:subB}
Let $V$ be a real vector space, let $E$ be a GC structure on $V$,
and let $E'$ be a $B$-field transform of $E$, defined by
$B\in\Wedge^2 V^*$. Then a subspace $W\subseteq V$ is a GC
subspace of $V$ with respect to $E$ if and only if it is a GC
subspace with respect to $E'$, and in that case, $E'_W$ is a
$B$-field transform of $E_W$ defined by $B\big\vert_W\in \Wedge^2
W^*$.
\end{cor}
\begin{proof}
Let $j:W\hookrightarrow V$ denote the inclusion map, and let
$\phi=\exp(u)\wedge f_1\wedge\dotsb\wedge f_k$ be a pure spinor
for the structure $E$, such that the $f_i$'s satisfy the
conditions of Proposition \ref{p:subspin}. It follows from
Proposition \ref{p:Bspin} that the spinor $\phi'=\exp(-B+u)\wedge
f_1\wedge\dotsb\wedge f_k$ is representative for $E'$. Now
Proposition \ref{p:subspin} implies that
\[
\phi_W:=\exp(j^* u)\wedge j^*(f_1)\wedge\dotsb\wedge j^*(f_l)
\quad \text{and}\quad \phi'_W:=\exp(-j^*B+j^* u)\wedge
j^*(f_1)\wedge\dotsb\wedge j^*(f_l)
\]
are representative spinors for $E_W$ and $E'_W$, respectively.
Since $B$ is real, it is immediate from Theorem
\ref{t:spinors}(\ref{i:pair}) that
$\pair{\phi_W}{\bar{\phi}_W}_M=\pair{\phi'_W}{\bar{\phi'}_W}_M$,
proving the first claim. The second assertion is clear.
\end{proof}

\subsection{GC quotients and duality}\label{ss:gcquot} Let $V$ be
a real vector space with a GCS defined by $E\subset V_{\bC}\oplus
V_{\bC}^*$. Consider a (real) subspace $W\subseteq V$ and the
corresponding quotient $V/W$. Let $\pi:V\to V/W$ be the projection
map, and let $\eta:\Ann(W)\to (V/W)^*$ be the natural isomorphism.
Define, dually to \eq{e:ew},
\begin{equation}\label{e:evw}
E_{V/W}=\Bigl\{ \bigl(\pi(\rho(e)),\eta(\rho^*(e))\bigr)
\,\big\vert\, e\in E\cap \bigl(V_{\bC}\oplus\Ann(W_{\bC})\bigr)
\Bigr\}.
\end{equation}
Again, it is easy to check (cf. Lemma \ref{l:dim}) that $E_{V/W}$
is a maximally isotropic subspace of
$(V/W)_{\bC}\oplus(V/W)_{\bC}^*$.
\begin{defin}\label{d:gcquot}
We say that $V/W$ is a {\em generalized complex quotient} of $V$
if $E_{V/W}\cap \overline{E}_{V/W}=(0)$. In this case, we denote
by $\cJ_{V/W}$ the orthogonal automorphism of $(V/W)\oplus(V/W)^*$
corresponding to $E_{V/W}$.
\end{defin}
Most questions about GC quotients can be easily reduced to
questions about GC subspaces by means of the following result.
\begin{prop}\label{p:subquot}
Let $V$ be a real vector space with a fixed GCS given by $E\subset
V_{\bC}\oplus V_{\bC}^*$. If $W\subseteq V$ is a real subspace,
then $V/W$ is a GC quotient of $V$ if and only if $\Ann(W)$ is a
GC subspace of $V^*$ with respect to $\tau(E)$ (cf.
\S\ref{ss:duality}). Suppose that this holds, and let $E_{V/W}$ be
the induced GCS on $V/W$. Let $E^D_{V/W}$ be the GCS on
$(V/W)^*\cong\Ann(W)$ induced by $\tau(E)\subset V_{\bC}^*\oplus
V_{\bC}$, and let $\tau_W:(V/W)\oplus (V/W)^*\to (V/W)^*\oplus
(V/W)$ be the isomorphism which interchanges the two summands.
Then $E^D_{V/W}=\tau_W(E_{V/W})$.
\end{prop}

\noindent The proposition follows immediately from Lemma
\ref{l:duality}. One might expect at first that $W$ is a GC
subspace of $V$ if and only if $V/W$ is a GC quotient of $V$. We
have discovered that this is not so: see Example
\ref{e:subnotquot}.

\subsection{Classification of GC vector spaces}\label{ss:classification}
We begin by defining direct sums of GC vector spaces. Let $U,V$ be
real vector spaces equipped with GC structures $E_U\subset
U_{\bC}\oplus U_{\bC}^*$, $E_V\subset V_{\bC}\oplus V_{\bC}^*$.
Let $\pi_U:U\oplus V\to U$, $\pi_V:U\oplus V\to V$ be the natural
projections, and let
\[
\nu_{U,V}:U\oplus U^*\oplus V\oplus V^* \rar{\cong} (U\oplus
V)\oplus (U\oplus V)^*
\]
denote the obvious isomorphism (or its complexified version). The
{\em direct sum} of the structures $E_U$ and $E_V$ is the GC
structure $\nu_{U,V}(E_U\oplus E_V)$ on $U\oplus V$. If $E_U$,
$E_V$ correspond to orthogonal automorphisms $\cJ_U$, $\cJ_V$,
then the automorphism corresponding to the direct sum is
$\nu_{U,V}\circ(\cJ_U\oplus\cJ_V)\circ\nu_{U,V}^{-1}$. Finally, if
$\phi_U$, $\phi_V$ are representative spinors for $E_U$ and $E_V$,
then $\pi_U^*(\phi_U)\wedge\pi_V^*(\phi_V)$ is a representative
spinor for the direct sum of the GC structures.

\sbr

The main result of the subsection is the following ``decomposition
theorem.''
\begin{thm}\label{t:classification}
Every GCS is a $B$-field transform of a direct sum of a complex
GCS and a symplectic GCS.
\end{thm}
This decomposition is by no means unique. Let us describe the
strategy of the proof of this theorem; along the way, we will give
a more canonical (but weaker) version of the result. Suppose that
$V$ is a GC vector space, with the GC structure defined by
$E\subset V_{\bC}\oplus V^*_{\bC}$. Let $S$ be the subspace of $V$
such that $S_{\bC}=\rho(E)\cap\rho(\overline{E})$. It is easy to
check that $S$ is a GC subspace of $V$, and in fact, the maximal
such subspace so that the induced structure on $S$ is
$B$-symplectic. This subspace is the ``$B$-symplectic part'' of
$V$, and is canonically determined. Moreover, $S$ doesn't change
if we make a $B$-field transform of $V$. The next step of the
proof is to choose an appropriate (non-unique) $B$-field transform
of the whole structure on $V$ such that $S$ becomes ``split'',
i.e., such that we can find a complementary GC subspace $W$ to $S$
in $V$ with $V=S\oplus W$ a direct sum of GC structures. To
complete the proof, we show that $W$ is $B$-complex.

\sbr

A dual construction is obtained by considering the subspace $C$ of
$V$ such that $C_{\bC}=(E\cap V_{\bC})\oplus (\overline{E}\cap
V_{\bC})$. It may not be a GC subspace of $V$ (see Example
\ref{ex:notquot}); on the other hand, $V/C$ is always a GC
quotient of $V$, and is $\be$-symplectic. Note also that $C$ is
obviously invariant under $\cJ$ (where $\cJ$ is the orthogonal
automorphism of $V\oplus V^*$ corresponding to $E$), and
$\cJ\restr{C}$ is a complex structure on $C$. Then we have the
following result.
\begin{thm}\label{t:canon}
\begin{enumerate}[(a)]
\item Every GC vector space $V$ has a unique maximal $B$-symplectic
subspace $S$. If $V/S$ is a GC quotient of $V$, then it is
$B$-complex. In any case, $V/S$ can be endowed with a canonical
complex structure.
\item Dually, every GC vector space $V$ has a unique minimal
subspace $C$ such that $V/C$ is a GC quotient which is
$\be$-symplectic. The subspace $C$ can be endowed with a canonical
complex structure, and $C$ satisfies the graph condition with
respect to this structure. If $C$ is a GC subspace, the induced
structure is $\be$-complex.
\end{enumerate}
\end{thm}
\begin{cor}\label{c:evendim}
If $V$ is a GC vector space, then $\dim_{\bR} V$ is even.
\end{cor}
\begin{proof}
By Theorem \ref{t:canon}(b), there exists a subspace $C\subseteq
V$ such that $C$ has a complex structure and $V/C$ is a
$\be$-symplectic GC quotient of $V$. Since a $\be$-symplectic
structure on a real vector space gives in particular a symplectic
structure, it follows that $\dim_{\bR} C$ and $\dim_{\bR}(V/C)$
are even. Hence, so is $\dim_{\bR} V$.
\end{proof}
We end by giving an explicit formula for a general $B$-field
transform of a direct sum of a symplectic and a complex GCS. Let
$S$ be a real vector space with a symplectic form $\om\in\Wedge^2
S^*$, let $C$ be a real vector space with a complex structure
$J\in\Aut_{\bR}(C)$, and form $V=S\oplus C$. We identify $V^*$
with $S^*\oplus C^*$ in the natural way. Let $B\in\Wedge^2 V^*$,
which we view as a skew-symmetric map $B:V\to V^*$, and write
accordingly as a matrix
\[
B=\matr{B_1}{B_2}{B_3}{B_4},
\]
where $B_1:S\to S^*$, $B_2:C\to S^*$, $B_3:S\to C^*$, $B_4:C\to
C^*$ are linear maps satisfying $B_1^*=-B_1$, $B_4^*=-B_4$,
$B_3=-B_2^*$. In the result that follows, an automorphism $\cJ$ of
$V\oplus V^*$ is viewed as a $4\times 4$ matrix according to the
decomposition $V\oplus V^*=S\oplus C\oplus S^*\oplus C^*$.
\begin{prop}\label{p:explicit}
With this notation, let $\cJ$ be the automorphism of $V\oplus V^*$
corresponding to the GCS on $V$ which is the $B$-field transform
of the direct sum of $(S,\om)$ and $(C,J)$ defined by $B$. Then
$\cJ$ is given by
\begin{equation}\label{e:Jsum}
\cJ=\left(
\begin{array}{cccc}
\om^{-1} B_1 & \om^{-1} B_2 & -\om^{-1} & 0 \\
0 & J & 0 & 0 \\
\om+B_1\om^{-1}B_1 & B_2 J + B_1\om^{-1}B_2 & -B_1\om^{-1} & 0 \\
B_3\om^{-1}B_1 + J^* B_3 & B_4 J + B_3\om^{-1} B_2 + J^* B_4 &
-B_3\om^{-1} & -J^*
\end{array}
\right).
\end{equation}
If $\pi_S:V\to S$, $\pi_C:V\to C$ are the two projections, then
the pure spinor corresponding to $\cJ$ is given by
\begin{equation}\label{e:phisum}
\phi=\exp(-B+i\pi_S^*\om)\wedge
\pi_C^*(f_1)\wedge\dotsb\wedge\pi_C^*(f_k),
\end{equation}
where $f_1,\dotsc,f_k$ are a basis of the $-i$-eigenspace of $J^*$
on $C^*_{\bC}$.
\end{prop}
The (entirely straightforward) proof is omitted. Observe that from
this description, it follows that $\om$, $J$ and $B_1$, $B_2$,
$B_3$ can be recovered from $\cJ$, while $B_4$ can only be
recovered up to the addition of a real bilinear skew form $\Om$ on
$C$ satisfying $\Om J+J^*\Om=0$.

\section{Split subspaces and the graph condition}\label{s:split}

\subsection{Summary} We begin the section by defining the twisting
of GC structures; we describe these both in terms of the
orthogonal automorphism $\cJ$ and in terms of representative
spinors. This notion is used to define the ``twisted product
structure'' on the direct sum of two GC vector spaces.

\sbr

Next we define ``generalized Lagrangian subspaces'' of GC vector
spaces; these were first introduced by Gualtieri \cite{Gua} and
Hitchin \cite{H3}, and were called ``generalized complex
subspaces'' by them. We also introduce a natural weakening of this
notion: ``generalized isotropic subspaces.'' These are used to
define another class of subspaces, namely, ``subspaces satisfying
the graph condition.''

\sbr

Finally, we introduce the most restricted class of subspaces: the
``split subspaces.'' We then study the relations between the
various types of subspaces that we have defined.

\subsection{Twisted GC structures}\label{ss:twist} Let $V$ be a real
vector space, and consider a GC structure on $V$ defined in terms
of an orthogonal automorphism $\cJ$ of $V\oplus V^*$ or the
corresponding pure spinor $\phi\in\Wedge\dot V_{\bC}^*$. As in
\S\ref{ss:notation}, we write $\cJ$ in matrix form:
\[
\cJ=\matr{\cJ_1}{\cJ_2}{\cJ_3}{\cJ_4}.
\]
We define the {\em twist} of $\cJ$ to be the GC structure on $V$
defined by
\[
\cJt=\matr{\cJ_1}{-\cJ_2}{-\cJ_3}{\cJ_4}.
\]
It is easy to check that $\cJt^2=-1$ and that $\cJt$ is orthogonal
with respect to the usual pairing $\pairing$ on $V\oplus V^*$. If
$V$ and $W$ are GC vector spaces, with the GC structure given by
$\cJ_V$ and $\cJ_W$, respectively, then by the {\em twisted
product structure} on $W\oplus V$ we will mean the direct sum (cf.
\S\ref{ss:classification}) of the GCS $\cJt_W$ on $W$ and the GCS
$\cJ_V$ on $V$.

\sbr

We can describe twisting in terms of spinors as follows.
\begin{prop}\label{p:twistspin}
Suppose that the GC structure $\cJ$ corresponds to the pure spinor
$\phi$, written in standard form as $\phi=\exp(u)\wedge
f_1\wedge\dotsb\wedge f_k$ (cf. Theorem \ref{t:spinors}). Then
$\cJt$ corresponds to the spinor $\tilde{\phi}=\exp(-u)\wedge
f_1\wedge\dotsb\wedge f_k$.
\end{prop}
\begin{proof}
Using our classification of GC structures (Theorem
\ref{t:classification}), we may assume that $\cJ$ has the form
\eqref{e:Jsum}. Then it it clear that twisting $\cJ$ is the same
as replacing $\om$, $B_1$, $B_2$, $B_3$ and $B_4$ by their
negatives. This operation replaces the spinor \eqref{e:phisum} by
the spinor $\exp(B-i\pi_S^*\om)\wedge
\pi_C^*(f_1)\wedge\dotsb\wedge\pi_C^*(f_k)$, whence the result.
\end{proof}

\begin{rem}
It is easy to check (without going into the proof of Proposition
\ref{p:twistspin}) that the operation
\[
\phi=\exp(u)\wedge f_1\wedge\dotsb\wedge f_k \mapsto
\tilde{\phi}=\exp(-u)\wedge f_1\wedge\dotsb\wedge f_k
\]
is well defined on pure spinors. Indeed, if
\[
\exp(u)\wedge f_1\wedge\dotsb\wedge f_k = \exp(v)\wedge
g_1\wedge\dotsb\wedge g_k,
\]
where $u,v$ are $2$-forms and $f_i$, $g_j$ are $1$-forms, then by
equating the homogeneous parts of degree $k+2m$ of both sides, we
find that
\[
u^m\wedge f_1\wedge\dotsb\wedge f_k = v^m\wedge
g_1\wedge\dotsb\wedge g_k
\]
for all $m\geq 0$, whence
\[
(-u)^m\wedge f_1\wedge\dotsb\wedge f_k = (-v)^m\wedge
g_1\wedge\dotsb\wedge g_k
\]
for all $m\geq 0$, proving our observation.
\end{rem}
\begin{examples}
\begin{enumerate}[(1)]
\item The twist of a complex GCS is the structure itself.
\item The twist of a symplectic GCS defined by a symplectic form $\om$
is the symplectic GCS defined by $-\om$.
\item If $\cJ$ is a GCS on a real vector space $V$ and $\cJ'$ is a
$B$-field transform of $\cJ$ defined by $B\in\Wedge^2 V^*$, then
$\cJt'$ is the $B$-field transform of $\cJt$ defined by $-B$.
\end{enumerate}
\end{examples}

\subsection{Generalized Lagrangian subspaces and the graph condition}\label{ss:lagrangian}
We begin by defining three types of subspaces of GCS, the last of
which was introduced by Gualtieri (\cite{Gua}) and Hitchin
(\cite{H3}).
\begin{defin}
Let $V$ be a real vector space with a GCS defined by
$\cJ\in\Aut_{\bR}(V\oplus V^*)$, and let $W\subseteq V$ be a
subspace.
\begin{enumerate}[(1)]
\item We say that $W$ is {\em generalized isotropic}
if $\cJ(W)\subseteq W\oplus\Ann(W)$.
\item We say that $W$ is {\em generalized coisotropic}
if $\cJ(\Ann(W))\subseteq W\oplus\Ann(W)$.
\item We say that $W$ is {\em generalized Lagrangian}
if $W$ is both generalized isotropic and generalized coisotropic,
that is, if $W\oplus\Ann(W)$ is stable under $\cJ$.
\end{enumerate}
\end{defin}

\begin{examples}
In the symplectic case, the three definitions specialize to usual
isotropic, coisotropic and Lagrangian subspaces, respectively. In
the complex case, all three definitions specialize to complex
subspaces (i.e., the subspaces stable under the automorphism $J$
defining the complex structure).\footnote{For more details, see
\S\ref{ss:exgenlag}.}
\end{examples}

\sbr

Recall from classical symplectic geometry (see, e.g., \cite{dS})
that if $(W,\om_W)$ and $(V,\om_V)$ are symplectic vector spaces,
then a vector space isomorphism $W\to V$ is a symplectomorphism if
and only if its graph is a Lagrangian subspace of $W\oplus V$ with
respect to the ``twisted product structure'' on the direct sum,
defined by the symplectic form $-\om_W+\om_V$. More generally,
suppose $(V,\om_V)$ is a symplectic vector space and $W\subseteq
V$ is a subspace equipped with some symplectic form $\om_W$. Then
$\om_W=\om_V\big\vert_W$ if and only if the graph of the inclusion
$W\hookrightarrow V$ is an isotropic subspace of $W\oplus V$ with
respect to the twisted product structure. This motivates the
following
\begin{defin}
Let $V$ be a GC vector space, and let $W\subseteq V$ be a subspace
equipped with a GC structure. We say that $W$ (together with the
given GCS) {\em satisfies the graph condition} if the graph of the
inclusion $W\hookrightarrow V$ is a generalized isotropic subspace
of $W\oplus V$ with respect to the twisted product structure.
\end{defin}

\subsection{Split subspaces}\label{ss:split} Let $V$ be a
real vector space with a GC structure defined by
$\cJ\in\Aut_{\bR}(V\oplus V^*)$. The following definition is
somewhat similar in spirit to the definition of generalized
Lagrangian subspaces.
\begin{defin}
We say that a subspace $W\subseteq V$ is {\em split} if there
exists a subspace $N\subseteq V$ such that $V=W\oplus N$ and
$W\oplus\Ann(N)$ is stable under $\cJ$.
\end{defin}
The terminology is explained by the first part of the following
result.
\begin{prop}\label{p:split}
Let $\cJ$ be a GCS on a real vector space $V$, and let $W\subseteq
V$ be a split subspace, so that $V=W\oplus N$ for some subspace
$N\subseteq V$ such that $W\oplus\Ann(N)$ is stable under $\cJ$.
Then
\begin{enumerate}[(a)]
\item Both $W$ and $N$ are GC subspaces of $V$. Moreover,
$V=W\oplus N$ is a direct sum of GC structures.
\item Consider the natural isomorphism $\psi:W\oplus\Ann(N)\to W\oplus W^*$,
given by $(w,f)\mapsto (w,f\big\vert_W)$. Then the induced GCS
$\cJ_W$ on $W$ has the form
$\cJ_W=\psi\circ\bigl(\cJ\big\vert_{W\oplus\Ann(N)}\bigr)\circ\psi^{-1}$.
\item The space $N\oplus\Ann(W)$ is also stable under $\cJ$,
and the induced structure $\cJ_N$ on $N$ has a similar
description.
\end{enumerate}
\end{prop}

\subsection{Relations between the various notions of subspaces}\label{ss:subspaces}
We have already seen six different notions of subspaces of GC
vector spaces: GC subspaces, subspaces satisfying the graph
condition, split subspaces, and generalized
isotropic/coisotropic/Lagrangian subspaces. We will now describe
in detail the relations between these types of subspaces. The
significance of generalized Lagrangian subspaces as graphs of
isomorphisms between GC vector spaces, and, more generally, the
role of generalized isotropic subspaces as the key ingredient in
our notion of the graph condition, has already been explained. The
rest is contained in the following three results.
\begin{prop}\label{p:splitgraph}
Suppose that $W$ is a split subspace of a GC vector space $V$.
Then $W$ satisfies the graph condition with the induced GC
structure.
\end{prop}
\begin{prop}\label{p:nascgraph}
Suppose that $V$ is a GC vector space with the GCS defined by
\[
\cJ=\matr{\cJ_1}{\cJ_2}{\cJ_3}{\cJ_4}.
\]
Then a GC subspace $W\subseteq V$ satisfies the graph condition
(with respect to the induced structure) if and only if $W$ is
stable under $\cJ_1$. Moreover, in this case, if $\cJ_W$ denotes
the induced GC structure on $W$, then
\[
\cJ_{W1}=\cJ_1\big\vert_W \text{  and  }
\cJ_{W3}(w)=\cJ_3(w)\big\vert_W\ \forall\, w\in W.
\]
\end{prop}
\begin{prop}\label{p:graphgcsub}
Suppose that $\cJ$ is a GCS on a real vector space $V$, and
suppose that $W\subseteq V$ is a subspace equipped with a GCS
$\cK$. Then
\begin{enumerate}
\item[(a)] $W$ satisfies the graph condition with respect to $\cK$ if
and only if $\cK_1=\cJ_1\restr{W}$ (in particular, $W$ is stable
under $\cJ_1$) and $\cK_3(w)=\cJ_3(w)\restr{W}$ for all $w\in W$.
\end{enumerate}
Now suppose that $W$ does satisfy the graph condition with respect
to $\cK$, and assume that $W$ is also a GC subspace\footnote{As we
show in Remark \ref{r:graphnotsub}, this condition is not
automatic.} of $V$. Then
\begin{enumerate}
\item[(b)] if $\cJ_W$ denotes the induced
GCS on $W$, then $W$ also satisfies the graph condition with
respect to $\cJ_W$; and
\item[(c)] $\cK$ is a $\be$-field transform of $\cJ_W$.
\end{enumerate}
\end{prop}
The last statement is based on the following
\begin{lem}\label{l:beta}
Let $\cJ$ and $\cJ'$ be two GC structures on a real vector space
$V$ which have the following forms:
\[
\cJ=\matr{\cJ_1}{\cJ_2}{\cJ_3}{\cJ_4} \text{  and  }
\cJ'=\matr{\cJ_1}{\cJ'_2}{\cJ_3}{\cJ_4}
\]
$($thus, only the $(1,2)$-entries are different$)$. Then there
exists a $\be\in\Wedge^2 V$ that transforms $\cJ$ into $\cJ'$.
Moreover, $\be$ is unique if we impose the additional constraint
that $\cJ_1\circ\be:V^*\to V$ be skew-symmetric $($here we view
$\be$ as a skew-symmetric linear map $V^*\to V$, as explained in
\S\ref{ss:notation}$)$.
\end{lem}

\section{Generalized complex manifolds and submanifolds}\label{s:gcm}

\subsection{Summary} We begin the section by recalling the
definitions of generalized almost complex and complex manifolds,
following \cite{Gua}, \cite{H2}, \cite{H3}. We slightly generalize
the approach of Hitchin and Gualtieri by discussing generalized
almost complex structures on vector bundles (other than the
tangent bundle to a manifold). We then explain how integrability
of a generalized complex structure on a manifold can be expressed
in terms of the {\em Courant bracket} or the {\em
Courant-Nijenhuis} tensor. We also mention that in this more
general approach, the right setting for defining integrability of
generalized almost complex structures seems to be the one provided
by the theory of Courant algebroids \cite{LWX}, \cite{Roy}.

\sbr

We then study the analogue of Courant's notion of an admissible
function \cite{Cou}, and we define the {\em Courant sheaf}, which
we show to be a sheaf of local Poisson algebras.

\sbr

Next we define the generalized complex submanifolds and split
submanifolds, as well as generalized isotropic (resp. coisotropic,
resp. Lagrangian) submanifolds of generalized complex manifolds.
The theory is completely analogous to its linear counterpart, in
that integrability here plays no role: the analogue of Courant's
theorem, stating that the Dirac structure on a submanifold induced
from an integrable Dirac structure on an ambient manifold is
itself integrable, also holds in our situation.

\sbr

We conclude the section by giving a justification for our notion
of a split submanifold. Namely, we show that the sections of a
Courant sheaf on a generalized complex manifold can be ``pulled
back'' to a split submanifold, in a natural (though non-canonical)
way. This was our original motivation for introducing the
definition of a split submanifold.

\subsection{Definitions of GAC manifolds}\label{ss:gacm}
Let $M$ be a real manifold. Following \cite{Gua} and
\cite{H1}--\cite{H3}, we define a {\em generalized almost complex
structure} (GACS) on a (finite rank) real vector bundle $V\to M$
to be one of the following two objects:
\begin{itemize}
\item A subbundle $E\subset V_{\bC}\oplus V_{\bC}^*$ which is
maximally isotropic with respect to the standard pairing
$\pairing$ and satisfies $E\cap\overline{E}=0$; or
\item An $\bR$-linear bundle automorphism $\cJ$ of $V\oplus V^*$
which is orthogonal with respect to $\pairing$ and satisfies
$\cJ^2=-1$.
\end{itemize}
The equivalence of these two descriptions is proved in the same
way as in the linear case (Proposition \ref{p:eqvdef}). If the
tangent bundle $TM$ of $M$ is endowed with a GACS, we call $M$ a
{\em generalized almost complex manifold} (GACM). It follows
immediately from Corollary \ref{c:evendim} that a GACM must have
even dimension as a real manifold.

\begin{rem}\label{r:gacspin}
One can also describe GAC structures on vector bundles in terms of
spinors. First observe that it is straightforward to generalize
the theory described in \S\ref{ss:spinors} to the case of vector
bundles on manifolds. Namely, if $V$ is a (finite rank) real
vector bundle on a (real) manifold $M$ and $V_{\bC}$ is its
complexification, one naturally defines the ``sheaf of Clifford
algebras'' $\cliff(V)$ which acts on the ``bundle of spinors''
$\Wedge\dot V_{\bC}^*$. If $E\subset V_{\bC}\oplus V_{\bC}^*$ is a
maximally isotropic subbundle, it is easy to check (using Theorem
\ref{t:spinors}(b)) that the subsheaf of $\Wedge\dot V_{\bC}^*$
consisting of the (germs of) sections that are annihilated by $E$
is a line subbundle of $\Wedge\dot V_{\bC}^*$. The Mukai pairing
is also defined and becomes a pairing $\pairing_M:\Wedge\dot
V_{\bC}^*\times\Wedge\dot V_{\bC}^*\to \Wedge^n V_{\bC}^*$, where
$n$ is the rank of $V$. Thus, alternatively, a GACS on a real
vector bundle $V\to M$ can be defined as the data of a line
subbundle $\Phi\subseteq\Wedge\dot V^*_{\bC}$ such that, for every
nonvanishing local section $\phi$ of $\Phi$, the function
$\pair{\phi}{\bar{\phi}}_M$ is nonvanishing, and $\phi$ is
pointwise a pure spinor in the sense of Definition \ref{d:puresp}.
\end{rem}

\subsection{Integrability}\label{ss:integrability} Recall that an
almost complex structure on a manifold $M$ arises from an actual
complex structure if a suitable integrability condition is
satisfied (this is the Newlander-Nirenberg theorem). Hitchin
defines ``generalized complex structures'' in terms of a natural
generalization of that integrability condition. To formulate this
condition, let us recall first the definition of the {\em Courant
bracket} (\cite{Cou}, p.~645) on the sections of $T_{\bC}M\oplus
T^*_{\bC}M$:
\[
\cou{X\oplus\xi}{Y\oplus\eta}=[X,Y]+\cL_X\eta-\cL_Y\xi+\frac{1}{2}\cdot
d(\iota_Y\xi-\iota_X\eta).
\]
Here, $[X,Y]$ is the usual Lie bracket on vector fields.
\begin{defin}[cf. \cite{Gua}, \cite{H2}, \cite{H3}]\label{d:gcm}
Let $M$ be a real manifold equipped with a GACS defined by
$E\subset T_{\bC}M\oplus T^*_{\bC}M$. We say that $E$ is {\em
integrable} if the sheaf of sections of $E$ is closed under the
Courant bracket. If that is the case, we also say that $E$ is a
{\em generalized complex structure} on $M$, and that $M$ is a {\em
generalized complex manifold} (GCM).
\end{defin}
\begin{rem}[cf. \cite{Cou},\cite{H3}]
The integrability condition for a GACS can be expressed in terms
of the corresponding automorphism $\cJ$, in a way completely
analogous to the usual almost complex case. Namely, integrability
of a GACS defined by $\cJ$ is equivalent to the vanishing of the
{\em Courant-Nijenhuis} tensor
\[
N_{\cJ}(X,Y)=\cou{\cJ X}{\cJ Y}-\cJ\cou{\cJ X}{Y}-\cJ\cou{X}{\cJ
Y}-\cou{X}{Y}
\]
where $X$, $Y$ are sections of $TM\oplus T^*M$.
\end{rem}
Integrability can also be expressed in terms of spinors; however,
we will not use this description in our paper. We refer the reader
to \cite{H3} for a discussion.
\begin{examples}[cf. \cite{Gua}]
\begin{enumerate}
\item Let $J$ be a usual almost complex structure on $M$. Applying
pointwise the construction of Definition \ref{d:complsympl}(a), we
get a generalized almost complex structure on $M$. This structure
is integrable in the sense of Definition \ref{d:gcm} if and only
if $J$ is integrable in the sense of the Newlander-Nirenberg
theorem.
\item Let $\om$ be a nondegenerate differential $2$-form on $M$. Applying
pointwise the construction of Definition \ref{d:complsympl}(c), we
get a generalized almost complex structure on $M$. This structure
is integrable if and only if $\om$ is closed, i.e., is a
symplectic form on $M$.
\end{enumerate}
\end{examples}

\noindent We also note that there are global analogues of the
notions of $B$- and $\be$-field transforms. For example, a {\em
closed} $2$-form $B$ acts on GC structures on $M$ in the same way
as described in \S\ref{ss:B-field}. We refer the reader to
\cite{Gua} for more details. Our terminology in the global case
will be the same as in the linear case; thus, for instance, we
will refer to a $B$-field transform of a complex (resp.,
symplectic) manifold as a {\em $B$-complex} (resp., {\em
$B$-symplectic}) GCM.

\begin{rem}\label{r:generalizeGCM}
Recall that in symplectic geometry, it is sometimes important to
consider symplectic forms on vector bundles over a manifold, other
than its tangent bundle (see \cite{Vais}). Thus, if $V$ is a
(finite rank) real vector bundle on a real manifold $M$, one says
that $V$ is a {\em symplectic vector bundle} if $V$ is equipped
with a nondegenerate $2$-form $\om\in\Ga(M,\Wedge^2 V^*)$. One can
then ask for a generalization of the closedness condition on
$\om$. An appropriate setting for this is provided by the notion
of a {\em Lie algebroid} (cf., e.g., \cite{Roy}). Recall that a
Lie algebroid structure on the vector bundle $V$ is given by a
bundle map $a:V\to TM$, called the {\em anchor}, and an
$\bR$-bilinear bracket $[\cdot,\cdot]_V$ on the sheaf of sections
of $V$, making the latter a sheaf of Lie algebras, and satisfying
the following conditions:
\begin{enumerate}[(1)]
\item $a\bigl([X,Y]_V\bigr)=[a(X),a(Y)]$ for all (local) sections
$X,Y$ of $V$, where $[\cdot,\cdot]$ is the usual Lie bracket on
vector fields;
\item $[X,fY]_V=f\cdot[X,Y]_V+(\cL_{a(X)}f)\cdot Y$ for all
sections $X,Y$ of $V$ and all (local) $C^\infty$ functions $f$ on
$M$.
\end{enumerate}
Such a structure gives a degree $1$ differential $d_V$ on the
sheaf of sections of $\Wedge\dot V^*$, making the latter a sheaf
of differential graded algebras. The definition of $d_V$ is given
by the usual Leibnitz formula. Explicitly (cf. \cite{Roy} or
\cite{NT}), if $\om$ is a local section of $\Wedge^p V^*$, then
\begin{eqnarray*}
(d_V\om)(X_0,\dotsc,X_p)&=&\sum_{i=0}^p (-1)^i \cL_{a(X_i)}
\om(X_0,\dotsc,\hat{X}_i,\dotsc,X_p) \\
&& +\sum_{0\leq i<j\leq p} (-1)^{i+j} \om([X,Y]_V,X_0,\dotsc,
\hat{X}_i,\dotsc,\hat{X}_j,\dotsc, X_p)
\end{eqnarray*}
for all appropriate sections $X_0,\dotsc,X_p$ of $V$.

\sbr

One can now define a {\em symplectic Lie algebroid} (cf.
\cite{NT}) to be a Lie algebroid $V$, together with a
nondegenerate $2$-form $\om\in\Ga(M,\Wedge^2 V^*)$, such that
$d_V\om=0$. Of course, by taking $V=TM$ with the usual bracket of
vector fields, we recover the usual notion of a symplectic
manifold.

\sbr

In view of these remarks, it seems that an appropriate general
setting for the study of generalized complex structures would be
the one provided by the theory of Courant algebroids (we refer the
reader to \cite{Roy} for a detailed study of these objects). We
will not give the definition of a Courant algebroid, as it is
rather long and is not used anywhere in our paper. Let us only
remark that, if $V$ and $V^*$ are a pair of Lie algebroids in
duality (\cite{LWX} or \cite{Roy}), then the direct sum $V\oplus
V^*$ acquires a natural structure of a Courant algebroid
(\cite{LWX}), and part of this structure is a bracket
$[\cdot,\cdot]$ on the sections of $V\oplus V^*$ which specializes
to the Courant bracket in case $V=TM$ (with the usual Lie bracket,
and with the anchor being the identity map) , $V^*=T^*M$ (with
zero bracket and anchor). Then one can define a GACS $E\subset
V_{\bC}\oplus V^*_{\bC}$ to be {\em integrable} if the sheaf of
sections of $E$ is closed under this bracket. In the case where
$V=TM$, one recovers Definition \ref{d:gcm}.
\end{rem}

\subsection{The Courant sheaf}\label{ss:coursheaf}
Let $M$ be a GACM, with the GACS defined by a maximally isotropic
subbundle $E\subset T_{\bC}M\oplus T^*_{\bC}M$. We define the {\em
Courant sheaf} of $M$ to be the following subsheaf $\Cou_M$ of
$C^\infty_M\oplus T_{\bC}M$:
\[
\Cou_M(U)=\bigl\{ (f,X)\in C^\infty_M(U)\oplus \Ga(U,T_{\bC}M)
\,\big\vert\, df+X\in\Ga(U,\overline{E}) \bigr\},
\]
for every open subset $U\subseteq M$. The origin of this
definition lies in Courant's notion (\cite{Cou}) of an admissible
function. The precise analogue of this notion in our situation is
the following: if $U\subseteq M$ is open, a function $f\in
C^\infty_M(U)$ is {\em admissible} if there exists an
$X\in\Ga(U,T_{\bC}M)$ with $(f,X)\in\Cou_M(U)$.

\sbr

The main properties of the Courant sheaf are summarized in the
following
\begin{prop}\label{p:cour}
Let $M$ be a generalized almost complex
manifold.
\begin{enumerate}[(a)]
\item The Courant sheaf is closed under componentwise addition,
and under the multiplication defined by
\[
(f,X)\cdot (g,Y) = (f\cdot g, f\cdot Y + g\cdot X).
\]
Thus, $\Cou_M$ is a sheaf of (commutative unital) $\bC$-algebras.
In particular, the product of two admissible functions is
admissible.
\item Furthermore, $\Cou_M$ is a sheaf of local rings. Namely, if
$m\in M$, then the maximal ideal of the stalk $\Cou_{M,m}$ is
formed by the germs of all sections $(f,X)$ of $\Cou_{M,m}$,
defined near $m$, such that $f(m)=0$.
\item Assume that $M$ is a generalized complex manifold. Then the
Courant sheaf is also closed under the bracket defined by
\[
\bigl\{ (f,X), (g,Y)\bigr\} = \Bigl( \frac{1}{2}\cdot
\bigl(X(g)-Y(f)\bigr), [X,Y]\Bigr) = \bigl( X(g), [X,Y] \bigr)
\]
$($the last equality uses the fact that $\overline{E}$ is
isotropic with respect to the standard pairing $\pairing${}$)$,
where $[\cdot,\cdot]$ is the usual Lie bracket of vector fields.
With this operation, $\Cou_M$ is a sheaf of Poisson algebras.
\end{enumerate}
\end{prop}
Part (c) of the proposition can be used to define a Poisson
bracket on the sheaf of admissible functions, as was done in
\cite{Cou}. In \S\ref{ss:splitsubman}, we will see that sections
of the Courant sheaf can be pulled back to a split submanifold.

\subsection{Submanifolds}\label{ss:subman} Let $S$ be a smooth
submanifold of a real manifold $M$, and let $E\subset
T_{\bC}M\oplus T_{\bC}^*M$ be a GACS. Applying pointwise the
construction of \S\ref{ss:gcsubspace} to the vector bundle
$TM\big\vert_S$ with the GACS $E\big\vert_S$ and the subbundle
$TS\subseteq TM\big\vert_S$, we obtain a {\em maximally isotropic
distribution} $E_S\subset T_{\bC}S\oplus T_{\bC}^*S$. That is,
$E_S$ is a subset of $T_{\bC}S\oplus T_{\bC}^*S$ such that the
fiber of $E_S$ over every point $s\in S$ is a maximally isotropic
subspace of $T_{s,\bC}M\oplus T^*_{s,\bC}M$.

\sbr

\noindent {\sc Caution}! In general, $E_S$ is {\em not} a
subbundle of $T_{\bC}S\oplus T_{\bC}^*S$.

\sbr

We refer the reader to \cite{Cou} for a detailed discussion of the
analogous situation for submanifolds of Dirac manifolds. In
particular, Courant's arguments can be easily adapted to give a
necessary condition under which $E_S$ is a subbundle of
$T_{\bC}S\oplus T_{\bC}^*S$ (cf. \cite{Cou}, Theorem 3.1.1), and
to prove that if $E_S$ is a subbundle and $E$ is integrable, then
so is $E_S$ (cf. \cite{Cou}, Corollary 3.1.4). This fact will not
be used in the sequel.
\begin{defin}\label{d:gcsubman}
We say that $S$ is a {\em generalized almost complex submanifold}
of $M$ if $E_S$ is a subbundle of $T_{\bC}S\oplus T_{\bC}^*S$ such
that $E_S\cap\overline{E}_S=0$. If moreover $M$ is a GCM, we will
simply say that $S$ is a {\em generalized complex submanifold} of
$M$.
\end{defin}

\sbr

We note, once again, that we have used different terminology than
the one introduced by Gualtieri and Hitchin in \cite{Gua},
\cite{H3}. They call a ``generalized complex submanifold'' what we
call a ``generalized Lagrangian submanifold'' in the definition
below.
\begin{defin}\label{d:gLsubman}
Let $M$ be a generalized (almost) complex manifold, with the
generalized (almost) complex structure defined by an orthogonal
automorphism $\cJ$ of $TM\oplus T^*M$. Let $S$ be a smooth
submanifold of $M$, and let $\Ann(TS)$ be the annihilator of the
subbundle $TS\subseteq TM\big\vert_S$ in $T^*M\big\vert_S$.
\begin{enumerate}[(1)]
\item We say that $S$ is a {\em generalized isotropic submanifold}
of $M$ if $\cJ(TS)\subseteq TS\oplus\Ann(TS)$.
\item We say that $S$ is a {\em generalized coisotropic
submanifold} of $M$ if $\cJ(\Ann(TS))\subseteq TS\oplus\Ann(TS)$.
\item We say that $S$ is a {\em generalized Lagrangian
submanifold} of $M$ if $S$ is both generalized isotropic and
generalized coisotropic, i.e., if $TS\oplus\Ann(TS)$ is stable
under $\cJ$.
\end{enumerate}
\end{defin}

\subsection{Split submanifolds}\label{ss:splitsubman}
Let $M$ be a generalized (almost) complex manifold, with the
generalized (almost) complex structure defined by an orthogonal
automorphism $\cJ$ of $TM\oplus T^*M$. We now define the global
analogue of the notion of a split subspace introduced in
\S\ref{ss:split}.
\begin{defin}
A smooth submanifold $S$ of $M$ is said to be {\em split} if there
exists a smooth subbundle $N$ of $TM\big\vert_S$ such that
$TM\big\vert_S=TS\oplus N$, and $TS\oplus\Ann(N)$ is invariant
under $\cJ$.
\end{defin}
The following result is the global version of, and immediately
follows from, Proposition \ref{p:split}.
\begin{prop}\label{p:splitsubman}
Let $\cJ\in\Aut_{\bR}(TM\oplus T^*M)$ be a GACS on a real manifold
$M$, and let $S\subset M$ be a split submanifold, so that
$TM\restr{S}=TS\oplus N$ for some subbundle $N\subseteq
TM\restr{S}$ such that $TS\oplus\Ann(N)$ is stable under $\cJ$.
Then $S$ is a GAC submanifold of $M$. Moreover, if
$\psi:TS\oplus\Ann(N)\to TS\oplus T^*S$ is the natural
isomorphism, then the induced GACS $\cJ_S$ on $S$ has the form
$\cJ_S=\psi\circ\bigl(\cJ\restr{TS\oplus\Ann(N)}\bigr)\circ\psi^{-1}$.
\end{prop}
The last result of the section relates split submanifolds to the
Courant sheaf:
\begin{thm}\label{t:pullcour}
Let $M$, $S$ and $N\subseteq TM\restr{S}$ be as above. Let
$\pi:T_{\bC}M\restr{S}\to T_{\bC}S$ denote the projection parallel
to $N_{\bC}$. For every open subset $U\subseteq M$, the map
\[
\Ga(U,C^\infty_M\oplus T_{\bC}M) \rar{} \Ga(U\cap S,
C^\infty_S\oplus T_{\bC}S),
\]
given by $(f,X)\mapsto \bigl(f\restr{S},\pi(X\restr{S})\bigr)$,
induces a ring homomorphism from $\Cou_M(U)$ to $\Cou_S(U\cap S)$.
In particular, the restriction of an admissible function to a
split submanifold is always admissible.
\end{thm}

\section{The generalized complex ``category''}\label{s:category}

\subsection{Summary} In this section we define another class of
submanifolds of generalized complex manifolds, namely,
``submanifolds satisfying the graph condition.'' This notion is
completely analogous to its linear counterpart studied in
\S\ref{s:split}. It relates our notion of GC submanifolds to the
notion of generalized isotropic subspaces, in the same way as
symplectic submanifolds of symplectic manifolds are related to
isotropic submanifolds (cf. \S\ref{ss:motivation}).

\sbr

We then use these ideas to define the generalized complex
``category''; it contains Weinstein's symplectic ``category''
\cite{Wei2} as a full ``subcategory.'' We end the section by
describing a simpler category -- the category of GC vector spaces
-- which is an honest category, and contains as a full subcategory
the category of symplectic vector spaces constructed by Guillemin
and Sternberg in \cite{GS}.

\subsection{The graph condition}\label{ss:submgraph}
We begin by defining twisted product structures for manifolds. Let
$M$, $N$ be real manifolds equipped with GAC structures defined by
$\cJ_M\in\Aut_{\bR}(TM\oplus T^*M)$, $\cJ_N\in\Aut_{\bR}(TN\oplus
T^*N)$. Applying pointwise the construction of \S\ref{ss:twist},
we obtain the notion of the {\em twist} of $\cJ_M$, which we will
denote by $\cJt_M$. It is easy to check that the twist of an
integrable GACS is again integrable. We then define the {\em
twisted product structure} on $M\times N$ to be the product of the
GAC structures defined by $\cJt_M$ and $\cJ_N$. (We leave to the
reader to define the notion of a product of GAC structures, which
is, again, completely analogous to the notion of a direct sum of
GC structures on vector spaces, defined in
\S\ref{ss:classification}. Note also that the product of two
integrable GAC structures is integrable.)

\sbr

Now let $M$ be a GAC manifold, and let $j:S\hookrightarrow M$ be a
smooth submanifold. Suppose that $S$ is equipped with some GACS.
We then say that $S$, together with this GACS, {\em satisfies the
graph condition}, if the graph of $j$ is a generalized isotropic
submanifold of $S\times M$ with respect to the twisted product
structure. We leave it to the reader to formulate the suitable
global analogues of the results of \S\ref{ss:subspaces}. Let us
only remark that this definition must be used with some caution
(cf. also \S\ref{ss:results} in the introduction). For instance,
if $j:S\to M$ is a diffeomorphism between two real manifolds, and
both $S$ and $M$ are equipped with GAC structures such that $j$
satisfies the graph condition, then it does not follow that $j$
induces an isomorphism between the two GAC structures. The reason
is that, in general (unlike the symplectic case) it is not true
that a generalized isotropic submanifold, which has half the
dimension of the ambient manifold, is generalized Lagrangian.
However, we have the following result, which motivates the
definition of the next subsection.
\begin{thm}\label{t:graphiso}
Let $\phi:M\to N$ be a diffeomorphism of real manifolds, and
suppose that $M$, $N$ are equipped with GAC structures. If the
graph of $\phi$ is generalized Lagrangian with respect to the
twisted product structure on $M\times N$, then $\phi$ induces an
isomorphism between the two GAC structures.
\end{thm}

\subsection{Definition of the GC ``category''}\label{ss:category}
We now introduce a ``category'' that generalizes Weinstein's
definition \cite{Wei2} of the symplectic ``category.'' Suppose
that $M$, $N$ are generalized complex manifolds. We define a {\em
canonical relation} from $M$ to $N$ to be a generalized Lagrangian
submanifold of $M\times N$ with respect to the twisted product
structure. Now suppose that $P$ is a third GCM, and let
$\Ga\subset M\times N$, $\Phi\subset N\times P$ be canonical
relations. The {\em set-theoretic composition} of $\Phi$ and $\Ga$
is the subset
\[
\Phi\circ\Ga = \bigl\{ (m,p)\in M\times P \,\big\vert\,
\text{there exists } n\in N \text{ with } (m,n)\in\Ga \text{ and }
(n,p)\in\Phi\bigr\}.
\]
If $\Phi\circ\Ga$ is a generalized Lagrangian submanifold of
$M\times P$ with respect to the twisted product structure, we say
that $\Phi$ and $\Ga$ are {\em composable}. This condition can be
expressed in terms of a certain transversality property. We then
define the {\em generalized complex ``category''} to be the
category whose objects are GC manifolds, and whose morphisms are
the canonical relations. The composition is only defined on pairs
of composable relations; thus, what we get is not a ``true''
category.

\sbr

The details of this construction will be discussed in our
subsequent publications.

\subsection{The category of GC vector spaces}\label{ss:catgcvs} We
finish the theoretical part of our paper by discussing the linear
counterpart of the category introduced in \S\ref{ss:category}. Our
exposition here follows \cite{GS}. If $V$, $W$ are real vector
spaces, a {\em relation} from $V$ to $W$ is a linear subspace of
$V\oplus W$. Let $Z$ be a third real vector space, and let
$\Ga\subset V\oplus W$, $\Phi\subset W\oplus Z$ be relations. The
{\em composition} of $\Phi$ and $\Ga$ is defined as is
\S\ref{ss:category} (see also \cite{GS}, p.~942):
\begin{equation}\label{e:comp}
\Phi\circ\Ga = \bigl\{ (v,z)\in V\oplus Z \,\big\vert\,
\text{there exists } w\in W \text{ with } (v,w)\in\Ga \text{ and }
(w,z)\in\Phi\bigr\}.
\end{equation}
It is clear that $\Phi\circ\Ga$ is a relation from $V$ to $Z$. The
main result of the subsection is the following
\begin{thm}\label{t:closcomp}
Let $V$, $W$, $Z$ be generalized complex vector spaces, and let
$\Ga\subset V\oplus W$, $\Phi\subset W\oplus Z$ be relations.
\begin{enumerate}[(a)]
\item If $\Ga$, $\Phi$ are generalized isotropic subspaces with
respect to the twisted product structures on $V\oplus W$ and
$W\oplus Z$, then $\Phi\circ\Ga$ is a generalized isotropic
subspace with respect to the twisted product structure on $V\oplus
Z$.
\item Same statement as part (a), with ``generalized isotropic'' replaced by
``generalized coisotropic.''
\item Same statement as part (a), with ``generalized isotropic'' replaced by
``generalized Lagrangian.''
\end{enumerate}
\end{thm}

This theorem allows us to introduce the following
\begin{defin}[The category $\Can$]\label{d:Can}
\begin{enumerate}[(1)]
\item If $V$, $W$ are GC vector spaces, a {\em canonical
relation} from $V$ to $W$ is a generalized Lagrangian subspace of
$V\oplus W$ with respect to the twisted product structure.
\item The {\em category $\Can$ of generalized complex vector spaces} is defined as
follows. The {\em objects} of $\Can$ are generalized complex
vector spaces. If $V,W$ are objects of $\Can$, the set of {\em
morphisms} from $V$ to $W$ is the set of canonical relations from
$V$ to $W$. The {\em composition} of morphisms is defined by
\eq{e:comp}.
\end{enumerate}
\end{defin}

\section{Examples and counterexamples}\label{s:examples}

\subsection{Summary}\label{ss:exsumm} In this section
we illustrate the main constructions and results of the previous
sections with examples (and counterexamples). Following the
philosophy of \cite{Gua}, we consider two basic types of examples:
$B$-complex and $B$-symplectic vector spaces. In particular, as a
result of our study, we will see that the three types of subspaces
we have introduced (GC subspaces, GC subspaces satisfying the
graph condition with respect to the induced structure, and split
subspaces) are in general all distinct.

\sbr

We begin by recalling explicit formulas for $B$- and $\be$-field
transforms of symplectic structures. We refer to \cite{Gua} for
the (straightforward) computations. Let $V$ be a real vector space
and let $B\in\Wedge^2 V^*$, $\be\in\Wedge^2 V$. If $J$ is a
complex structure on $V$ and $\cJ_J$ is the corresponding GCS (see
Definition \ref{d:complsympl}(a)), then the transforms of $\cJ$
defined by $B$ and $\be$ are given by
\[
\matr{J}{0}{BJ+J^*B}{-J^*} \ \ \ \text{  and  }\ \ \
\matr{J}{-J\be-\be J^*}{0}{-J^*},
\]
respectively. If $\om$ is a symplectic form on $V$ and $\cJ_{\om}$
is the corresponding GCS (see Definition \ref{d:complsympl}(c)),
then the transforms of $\cJ$ defined by $B$ and $\be$ are given by
\[
\matr{\om^{-1}B}{-\om^{-1}}{\om+B\om^{-1}B}{-B\om^{-1}} \ \ \
\text{  and  }\ \ \
\matr{\be\om}{-\be\om\be-\om^{-1}}{\om}{-\om\be},
\]
respectively.

\sbr

It is also easy to describe $B$-field transforms of complex and
symplectic structures in terms of the subspace $E\subset
V_{\bC}\oplus V_{\bC}^*$. Let us start with the complex case. If
$J$ is a complex structure on $V$, we call the elements of the
$+i$-eigenspace (resp., $-i$-eigenspace) of $J$ in $V_{\bC}$ the
{\em holomorphic} (resp., {\em antiholomorphic}) vectors; and we
call the elements of the $+i$-eigenspace (resp., $-i$-eigenspace)
of $J^*$ in $V_{\bC}^*$ the {\em holomorphic} (resp., {\em
antiholomorphic}) $1$-forms on $V_{\bC}$. Then the $B$-field
transform of the complex GCS defined by $J$ is given by
\[
E=\bigl\{ (v,\iota_v B+ f) \,\big\vert\, v\in V_{\bC} \text{ is
holomorphic, } f\in V^*_{\bC} \text{ is antiholomorphic}\bigr\}.
\]
On the other hand, if $\om$ is a symplectic form on $V$, then the
$B$-field transform of the symplectic GCS defined by $\om$ is
given by
\[
E=\bigl\{ (v, \iota_v(B-i\om)) \,\big\vert\, v\in V_{\bC} \bigr\}.
\]

\subsection{Subspaces, quotients and submanifolds}\label{ss:exsubquot}
Let us classify GC subspaces of $B$-complex and $B$-symplectic
vector spaces.
\begin{prop}\label{p:subBcs}
Let $V$ be a real vector space, let $B\in\Wedge^2 V^*$, and let
$W\subseteq V$ be a subspace.
\begin{enumerate}[(a)]
\item If $J$ is a complex structure on $V$ and $\cJ_J$ is the
associated GCS, then $W$ is a GC subspace of $V$ with respect to
the transform of $\cJ_J$ defined by $B$ if and only if $W$ is
invariant under $J$. In this case, the induced structure on $W$ is
also $B$-complex, namely, it is the transform by $B\big\vert_W$ of
the complex GCS associated to $J\big\vert_W$. The induced
structure satisfies the graph condition. Furthermore, $W$ is split
if and only if there exists a complement $N$ to $W$ in $V$ such
that $N$ is also invariant under $J$, and $N$ is orthogonal to $W$
with respect to $B$.
\item If $W$ is a subspace of a $B$-complex GC vector space which
satisfies the graph condition with respect to some GC structure,
then $W$ is also a GC subspace.
\item If $\om$ is a symplectic form on $V$ and $\cJ_\om$ is the
associated GCS, then $W$ is a GC subspace of $V$ with respect to
the transform of $\cJ_\om$ defined by $B$ if and only if
$\om\restr{W}$ is nondegenerate. In this case, the induced
structure on $W$ is also $B$-symplectic, namely, it is the
transform by $B\big\vert_W$ of the symplectic GCS associated to
$\om\big\vert_W$. Furthermore, $W$ is split if and only if the
induced structure satisfies the graph condition, and either of
these is equivalent to: if $N$ is the orthogonal complement to $W$
in $V$ with respect to $\om$, then $N$ is also orthogonal to $W$
with respect to $B$.
\end{enumerate}
\end{prop}
The analogue of part (b) holds in the symplectic case, but fails
in the more general $B$-symplectic case (cf. Remark
\ref{r:graphnotsub}). Using duality, it is easy to obtain a
similar description of quotients of $\be$-complex and
$\be$-symplectic structures. We omit the simple proof of
Proposition \ref{p:subBcs}, which is a combination of the results
of Propositions \ref{p:splitgraph}--\ref{p:graphgcsub} and
Corollary \ref{c:subB}, and the explicit formulas of
\S\ref{ss:exsumm}.

\sbr

This classification of GC subspaces generalizes without difficulty
to a classification of GC submanifolds of $B$-complex and
$B$-symplectic manifolds. Thus, for example, if $M$ is endowed
with a $B$-field transform of a complex structure
$J\in\Aut_{\bR}(TM)$, then a submanifold $j:S\hookrightarrow M$ is
generalized complex if and only if $TS\subseteq TM\big\vert_S$ is
stable under $J$. It is easy to check that the condition that
$E_S$ be a smooth subbundle of $T_{\bC}S\oplus T_{\bC}^*S$ (cf.
\S\ref{ss:subman}) is automatic. The induced structure on $S$ is
also $B$-complex, defined by $J\big\vert_{TS}$ and $j^* B$. A
similar remark applies in the $B$-symplectic case. We leave to the
reader the task of formulating the suitable descriptions of
submanifolds satisfying the graph condition, and of split
submanifolds, completely analogous to Proposition \ref{p:subBcs}.

\sbr

Finally, we observe that all GC submanifolds of complex or
symplectic GC manifolds are already split. In the symplectic case,
the ``splitting'' is canonical. More precisely, if $S$ is a GC
submanifold of a symplectic manifold $(M,\om)$, then there exists
a unique subbundle $N\subset TM\big\vert_S$ such that
$TM\big\vert_S=TS\oplus N$ is a direct sum of GC structures;
namely, $N$ is the orthogonal complement to $TS$ with respect to
$\om$. In the complex case, the ``splitting'' may not be unique,
but there always exists one. Indeed, let $J$ be a complex
structure on a real manifold $M$, and let $S$ be a submanifold of
$M$ such that $TS\subseteq TM\big\vert_S$ is $J$-invariant. Then
$J$ also acts on the normal bundle $\cN=(TM\big\vert_S)/TS$, so we
get a short exact sequence
\[
0 \rar{} TS \rar{} TM\restr{S} \rar{} \cN \rar{} 0.
\]
Since we are working in the $C^\infty$ category, there is a smooth
splitting of this sequence $\sg:\cN\to TM\restr{S}$. We only need
to make $\sg$ equivariant with respect to $J$; since $J^4=1$, this
is easily achieved by replacing $\sg$ with $(1/4)\cdot
\sum_{k=0}^3 J^k\circ\sg\circ J^{-k}$.

\subsection{Generalized Lagrangian subspaces and submanifolds}\label{ss:exgenlag}
We now briefly discuss generalized isotropic, coisotropic, and
Lagrangian subspaces (resp., submanifolds) of $B$-complex and
$B$-symplectic GC vector spaces (resp., GC manifolds). First let
$J$ be a complex structure on a real vector space $V$, inducing a
complex GCS as usual. From Definition \ref{d:complsympl}(a), it is
clear that a subspace $W\subseteq V$ is generalized isotropic if
and only if $J(W)\subseteq W$, and it is generalized coisotropic
if and only if $J^*(\Ann(W))\subseteq\Ann(W)$. The latter
condition is again equivalent to $J(W)\subseteq W$. Thus, in the
complex case, the following classes of subspaces are all the same:
generalized isotropic, generalized coisotropic, generalized
Lagrangian, and generalized complex; and they all coincide with
the class of subspaces invariant under $J$.

\sbr

Next let $\om$ be a symplectic form on $V$, inducing a symplectic
GCS as usual. From Definition \ref{d:complsympl}(c), it is clear
that a subspace $W\subseteq V$ is generalized isotropic if and
only if $\om(W)\subseteq\Ann(W)$, i.e., $\om\big\vert_W\equiv 0$,
that is, $W$ is isotropic in the sense of classical symplectic
geometry. Similarly, $W$ is generalized coisotropic if and only if
$W$ is coisotropic in the classical sense; and hence $W$ is
generalized Lagrangian if and only if $W$ is Lagrangian.

\sbr

These remarks carry over to the global case without any changes.
Thus, if $M$ is a complex GCM, then the generalized
isotropic/coisotropic/Lagrangian submanifolds of $M$ are precisely
the complex submanifolds in the usual sense. Similarly, if $M$ is
a symplectic GCM, then the generalized
isotropic/coisotropic/Lagrangian submanifolds of $M$ are the
submanifolds that are isotropic/coisotropic/Lagrangian in the
classical sense.

\sbr

Finally, we consider the case of $B$-complex and $B$-symplectic GC
manifolds. For simplicity, we restrict ourselves to generalized
Lagrangian submanifolds, which are studied in detail (under the
name of ``generalized complex submanifolds'') in \cite{Gua}. Thus,
we only list the results, following \cite{Gua}.
\begin{itemize}
\item If $M$ is endowed with a GCS which is a transform of a
complex structure on $M$ by a $B$-field $B\in\Ga(M,\Wedge^2
T^*M)$, then a submanifold $j:S\hookrightarrow M$ is a generalized
Lagrangian submanifold if and only if $S$ is a complex submanifold
of $M$, and $j^* B$ is a form of type $(1,1)$ on $S$.
\item If $M$ is endowed with a GCS which is a transform of a
symplectic structure on $M$ by a $B$-field $B\in\Ga(M,\Wedge^2
T^*M)$, then a submanifold $j:S\hookrightarrow M$ is a generalized
Lagrangian submanifold if and only if $S$ is a coisotropic
submanifold of $M$ (in the classical sense), and, if $TS^\perp$
denotes the orthogonal complement to $TS$ in $TM\big\vert_S$ with
respect to the symplectic form on $M$ , then $TS$ is orthogonal to
$TS^\perp$ with respect to $B$, so that $j^* B$ descends to
$TS/TS^\perp$; and, finally, $(j^*\om)^{-1}j^* B$ is an almost
complex structure on $TS/TS^\perp$.
\end{itemize}

\subsection{Symplectic geometry with $B$-fields}\label{ss:symplB}
Let $V$ be a real vector space, and let $\om$ be a symplectic form
on $V$. We have an injective linear map
\[
\Wedge^2 V^* \rar{} \End_{\bR}(V),
\]
given by $B\mapsto T:=\om^{-1}B$. A necessary and sufficient
condition for $T\in\End_{\bR}(V)$ to be in the image of this map
is that $(\om T)^*=-\om T$, or, equivalently,
\begin{equation}\label{e:star}
\om(u,Tv)=\om(Tu,v)\ \ \forall\,u,v\in V.
\end{equation}
If $T$ satisfies \eq{e:star}, so that $T=\om^{-1}B$ for some
$B\in\Wedge^2 V^*$, then the $B$-field transform of the symplectic
GCS defined by $\om$ can be written in terms of $T$ as follows:
\[
\matr{T}{-\om^{-1}}{\om\circ(1+T^2)}{-T^*}.
\]
Thus, from Proposition \ref{p:characterization}, we immediately
get the following
\begin{prop}
The $B$-field transform of $\om$ is
\begin{itemize}
\item symplectic $\iff$ $T=0$;
\item $\be$-symplectic $\iff$ $i$ is not an eigenvalue of $T$;
\item $\be$-complex $\iff$ $T^2=-1$.
\end{itemize}
\end{prop}
The second possibility can be realized, e.g., by taking $T=1$
(which corresponds to $B=\om$). This shows that a (nontrivial)
$B$-field transform of a symplectic GCS may end up being
$\be$-symplectic as well. To show that the third possibility can
also be realized, let us choose a basis
$\{e_1,\dotsc,e_n,f_1,\dotsc,f_n\}$ of $V$ in which $\om$ has the
canonical form:
\begin{equation}\label{e:symplcan}
\om(e_i,f_j)=\de_{ij},\ \ \om(e_i,e_j)=\om(f_i,f_j)=0\
\forall\,i,j.
\end{equation}
Pick any $T\in\End_{\bR}(V)$, and write $T$ with respect to this
basis in block form as follows:
\begin{equation}\label{e:T}
T=\matr{A}{B}{C}{D}.
\end{equation}
Then \eq{e:star} is equivalent to
\[
\matr{0}{I}{-I}{0}\cdot\matr{A}{B}{C}{D} =
\matr{A^t}{C^t}{B^t}{D^t}\cdot\matr{0}{I}{-I}{0},
\]
where $I$ is the $n\times n$ identity matrix, and superscript $t$
denotes the transpose. The last equality can be rewritten as
\begin{equation}\label{e:condT}
D=A^t,\ \ B^t=-B,\ \ C^t=-C.
\end{equation}
In particular, if $n$ is even, we can choose $A$ so that $A^2=-1$;
taking $B=C=0$, we get an example of an operator $T$ satisfying
\eq{e:star} such that $T^2=-1$.

\sbr

\noindent {\sc Observation.} {\em This provides an example of a
situation where a $B$-field transform of a symplectic GCS is also
$\be$-complex.}

\begin{rem}\label{r:graphnotsub}
As a by-product of our discussion, we also obtain an example of a
subspace of a $B$-symplectic vector space which satisfies the
graph condition with respect to some GC structure, but is not a GC
subspace in the sense of Definition \ref{d:gcsub}. For example, if
$W$ is the subspace of $V$ spanned by $e_1,\dotsc,e_n$, then $W$
is invariant under the operator $T$ defined above, so $T\restr{W}$
gives a complex structure on $W$. It is then easy to see (using
Proposition \ref{p:graphgcsub}) that $W$ satisfies the graph
condition with respect to this complex structure. However, $W$ is
obviously not a GC subspace of $V$ (cf. Proposition
\ref{p:subBcs}), since $\om\restr{W}\equiv 0$.
\end{rem}

\sbr

Now let us pass to the global case. Thus $\om$ is a symplectic
form on a real manifold $M$, and $T$ is an $\bR$-linear bundle
endomorphism of $TM$ satisfying \eq{e:star} for all sections $u,v$
of $TM$. We suppose that if $B$ is the two-form on $M$
corresponding to $T$, then the $B$-field transform of the GCS
defined by $\om$ is integrable. Let $E\subset T_{\bC}M\oplus
T_{\bC}^*M$ denote the maximally isotropic subbundle corresponding
to this transform. It is clear that $E\cap T_{\bC}M$ is the
$+i$-eigensheaf of $T$ on $T_{\bC}M$. Clearly this sheaf is closed
under the Lie bracket of vector fields (because the restriction of
the Courant bracket to $T_{\bC}M$ is the usual Lie bracket). The
same applies to $\overline{E}\cap T_{\bC}M$. Now $(E\cap
T_{\bC}M)\oplus(\overline{E}\cap T_{\bC}M)$ is the
complexification of the subsheaf $\Ker(1+T^2)\subseteq TM$, whence
the latter subsheaf is integrable.

\sbr

There exists a dense open set $U\subseteq M$ on which $1+T^2$ has
constant rank. On this set, the three sheaves considered above
become subbundles of $T_{\bC}M$. Thus we obtain
\begin{cor}
There exists a natural foliation of $U$, corresponding to the
integrable smooth distribution $\Ker(1+T^2)$, whose leaves acquire
a natural complex structure.
\end{cor}
This corollary is a global analogue of a special case of Theorem
\ref{t:canon}(b).

\subsection{Some counterexamples}\label{ss:counterex}
We start by giving an example of a $B$-symplectic vector space $V$
and a GC subspace $W\subset V$ such that $V/W$ is not a GC
quotient of $V$.
\begin{example}\label{e:subnotquot}
Let $V$ be a real vector space with basis $p_1,q_1,p_2,q_2$, and
let $a_1,b_1,a_2,b_2$ be the dual basis of $V^*$. Define
$\om=a_1\wedge b_1+a_2\wedge b_2$, a symplectic form on $V$, and
$B=a_1\wedge a_2 - b_1\wedge b_2$, a $B$-field on $V$. Let $W$ be
the plane in $V$ spanned by $p_1$ and $q_1$. Then $W$ is a GC
subspace of $V$, but $V/W$ is not a GC quotient of $V$.
\end{example}
Indeed, it is clear that $W$ is a GC subspace of $V$, because
$\om\big\vert_W$ is nondegenerate. Let $\pi:V\to V/W$ denote the
quotient map, and let $\eta:\Ann(W)\to (V/W)^*$ be the natural
isomorphism. It follows from \S\ref{ss:exsumm} and \eqref{e:evw}
that
\[
E_{V/W}=\Bigl\{ \bigl(\pi(v), \eta(\iota_v(B-i\om))\bigr)
\,\Big\vert\, \iota_v(B-i\om)\in\Ann(W_{\bC}) \Bigr\}.
\]
Consider the element $v=p_1+i q_2$. It is easy to compute that
$\iota_v(B-i\om)=0$. Thus, $(\pi(v),0)\in E_{V/W}$. But since
$p_1\in W$, we have
$\overline{\pi(v)}=\overline{\pi(iq_2)}=\pi(-iq_2)=-\pi(v)$,
whence we also have $(\pi(v),0)\in\overline{E}_{V/W}$. This
implies that $E_{V/W}\cap\overline{E}_{V/W}\neq (0)$, so $V/W$ is
not a GC quotient of $V$.

\sbr

Now we discuss a more complicated example.
\begin{example}\label{ex:notquot}
There exists a $B$-symplectic vector space $V$ with the following
property. If $C$ is the minimal subspace of $V$ such that $V/C$ is
a $\be$-symplectic quotient of $V$ (cf. Theorem \ref{t:canon}),
then $C$ is not a GC subspace of $V$. Hence, by duality, we see
that there exists a $\be$-symplectic vector space $V$ such that,
if $S$ is the maximal $B$-symplectic subspace of $V$, then $V/S$
is not a GC quotient of $V$.
\end{example}
To find an example of such a $V$, we use the theory developed in
\S\ref{ss:symplB}. We suppose that $V$ has a distinguished basis
$\{e_i,f_i\}_{i=1}^n$, and is equipped with a symplectic form
$\om$ defined by \eqref{e:symplcan}. We consider an operator $T$
on $V$, represented by the matrix \eqref{e:T} (where $A$, $B$,
$C$, $D$ are to be chosen later). It is clear that with this
notation, the subspace $C\subset V$ is nothing but the kernel of
$1+T^2$. Now we note that the equation \eq{e:star} implies that
$\Ker(1+T^2)$ is orthogonal to $\Image(1+T^2)$ with respect to
$\om$. Thus, if $\Ker(1+T^2)\cap\Image(1+T^2)\neq (0)$, the
restriction of $\om$ to $\Ker(1+T^2)$ is degenerate, so that
$C=\Ker(1+T^2)$ is not a GC subspace of $V$.

\sbr

It remains to choose an arbitrary $n\times n$ matrix $A$ such that
$\Ker(1+A^2)\cap\Image(1+A^2)\neq (0)$; for instance, one can take
\[
A=\matr{J}{I}{0}{J},
\]
where
\[
J=\matr{0}{1}{-1}{0}, \quad I=\matr{1}{0}{0}{1},
\]
so that
\[
1+A^2=\matr{0}{2J}{0}{0}
\]
is a nonzero nilpotent matrix. Setting $D=A^t$, $B=C=0$, we obtain
an example of an operator $T$ of the form \eqref{e:T}, satisfying
\eqref{e:condT}, and such that $\Ker(1+T^2)\cap\Image(1+T^2)\neq
(0)$.
\begin{rem}
We leave it to the reader to check that an example of the type
presented above does not exist in the $B$-complex case. That is,
if $V$ is a $B$-complex GC vector space, and $C$ is the minimal
subspace such that $V/C$ is a $\be$-symplectic GC quotient of $V$,
then $C$ is a GC subspace of $V$, and the induced structure on $C$
is complex.
\end{rem}

\section{Proofs}\label{s:proofs}

\subsection{Summary} In this section, we give proofs of the
results that were stated (and not proved) in
\S\S\ref{s:linalg}--\ref{s:category}. This being a preliminary
version of the paper, in many places we only give sketches of the
arguments.

\subsection{Proofs for section
\ref{s:linalg}}\label{ss:proofslinalg} We start by proving
Proposition \ref{p:eqvdef}. A bijection between the collection of
$\cJ\in\Aut_{\bR}(V\oplus V^*)$ and the collection of $E\subset
V_{\bC}\oplus V_{\bC}^*$ is obtained by associating to $\cJ$ its
$+i$-eigenspace in $V_{\bC}\oplus V^*_{\bC}$. [Conversely, if we
start with $E$, we consider the $\bC$-linear automorphism of
$V_{\bC}\oplus V_{\bC}^*$ defined as multiplication by $+i$ on $E$
and multiplication by $-i$ on $\overline{E}$; it is
straightforward to check that this automorphism is induced by an
$\bR$-linear automorphism of $V\oplus V^*$.] On the other hand, a
bijection between the collection of pure spinors
$\phi\in\Wedge\dot V_{\bC}^*$ and the collection of $E$ is
obtained by associating to $\phi$ its annihilator in
$V_{\bC}\oplus V^*_{\bC}$ with respect to the Clifford action (cf.
\S\ref{ss:spinors}). The fact that this is a bijection follows
from parts (\ref{i:ann}) and (\ref{i:conj}) of Theorem
\ref{t:spinors}.

\sbr

For the subsequent arguments, we will need an explicit description
of the inverse of the bijection between $\phi$'s and $E$'s.
Suppose that $\phi\in\Wedge\dot V^*_{\bC}$ is a pure spinor,
written in the form of Theorem \ref{t:spinors}(\ref{i:pure}):
$\phi=\exp(u)\wedge f_1\wedge\dotsb\wedge f_k$. Let $\Phi$ denote
the span of $f_1,\dotsc,f_k$; this is a $k$-dimensional complex
subspace of $V^*_{\bC}$. We claim that
\begin{equation}\label{e:ann}
E_\phi:=\bigl\{ v-\iota_v(u)+f \,\big\vert\, v\in\Ann(\Phi),\
f\in\Phi \bigr\}
\end{equation}
is the maximally isotropic subspace corresponding to $\phi$.
Indeed, it is obvious that $\dimc E_\phi=\dimc V_{\bC}$, so it
suffices to see that all elements of $E$ annihilate $\phi$. By
definition of the Clifford action, we have
$f\cdot\phi=f\wedge\phi=0$ for all $f\in\Phi$. Next, if
$v\in\Ann(\Phi)$, then, using the fact that $\iota_v$ is a
superderivation of the algebra $\Wedge\dot V_{\bC}^*$, we find
that
\[
\iota_v(\phi) = \iota_v(\exp(u))\wedge f_1\wedge\dotsb\wedge f_k =
\iota_v(u)\wedge\exp(u)\wedge f_1\wedge\dotsb\wedge f_k =
\iota_v(u)\wedge\phi,
\]
whence also
$(v-\iota_v(u))\cdot\phi=\iota_v(\phi)-\iota_v(u)\wedge\phi=0$. We
can now give the

\begin{proof}[Proof of Proposition \ref{p:Bspin}]
Let $\cB$ be the orthogonal automorphism of $V\oplus V^*$
corresponding to $B$, as defined in \S\ref{ss:B-field}.
Explicitly, $\cB$ acts as follows: $\cB(v+f)=v+\iota_v(B)+f$. Now
let us write $\phi$ in the canonical form: $\phi=\exp(u)\wedge
f_1\wedge\dotsb\wedge f_k$. It is clear from Theorem
\ref{t:spinors}(\ref{i:pure}) that $\psi:=\exp(-B)\wedge\phi$ is
also a pure spinor. Using \eq{e:ann}, we immediately find that
$E_\psi=\cB\cdot E_\phi$, which proves the proposition.
\end{proof}

Before we proceed, we list some useful formulas. If $V$ is a real
vector space and $\cJ$ is an $\bR$-linear automorphism of $V\oplus
V^*$, written, as usual, in matrix form:
\[
\cJ=\matr{\cJ_1}{\cJ_2}{\cJ_3}{\cJ_4},
\]
then the conditions that $\cJ$ preserves the pairing $\pairing$
and satisfies $\cJ^2=-1$ are equivalent to the following
collection of equations:
\begin{eqnarray}
\label{e:1} &&\cJ_1^2+\cJ_2\cJ_3=-1; \\
\label{e:2} &&\cJ_1\cJ_2+\cJ_2\cJ_4=0; \\
\label{e:3} &&\cJ_3\cJ_1+\cJ_4\cJ_3=0; \\
\label{e:4} &&\cJ_4^2+\cJ_3\cJ_2=-1; \\
\label{e:5} &&\cJ_4=-\cJ_1^*; \\
\label{e:6} &&\cJ_2^*=-\cJ_2; \\
\label{e:7} &&\cJ_3^*=-\cJ_3.
\end{eqnarray}
The verification of this claim is straightforward. It is used in
the proof below.
\begin{proof}[Proof of Proposition \ref{p:characterization}]
We first note that in parts (a) and (c), it is enough to prove the
equivalence of the given statements involving $B$-field
transforms, because then the analogous result for $\be$-field
transforms follows formally with the aid of duality
(\S\ref{ss:duality} and Remark \ref{r:dualityBbeta}). Moreover,
part (b) is immediate from the explicit formulas of
\S\ref{ss:exsumm}. For part (d), note that if $\cJ_1=0$, then the
equations \ref{e:1}--\ref{e:7} imply that $\cJ_4=0$ and
$\cJ_2=-\cJ_3^{-1}$, with $\cJ_2^*=-\cJ_2$. Hence, taking
$\om=\cJ_3$, we see from Definition \ref{d:complsympl} that $E$ is
symplectic.

\sbr

The equivalences
\[
\rho(E)\cap\rho(\overline{E})=(0)\ \iff\ V_{\bC}^*=(V_{\bC}^*\cap
E)+(V_{\bC}^*\cap\overline{E})\ \iff\ \cJ_2=0
\]
and
\[
E\cap V_{\bC}^*=(0)\ \iff\ \rho(E)=V_{\bC}\ \iff\ \cJ_2 \text{ is
an isomorphism}
\]
are immediate.

\sbr

Now if $E$ is $B$-complex (resp., $B$-symplectic), it is clear
from \S\ref{ss:exsumm} that $\cJ_2=0$ (resp., $\cJ_2$ is an
isomorphism). Conversely, suppose that $\cJ_2=0$. Then it is easy
to check, using the formulas \eq{e:1}--\eq{e:7}, that $\cJ_1^2=-1$
and $\cJ_4=-\cJ_1^*$. We put $J=\cJ_1$; thus $J$ is a complex
structure on $V$. If we let $B=-(1/2)\cdot\cJ_3\cdot\cJ_1$, it is
straightforward to check that $B$ is skew-symmetric and
$\cJ_3=B\cJ_1+\cJ_1^*B$. Thus $E$ is $B$-complex.

\sbr

Finally, suppose that $\cJ_2$ is an isomorphism. We put
$\om=-\cJ_2^{-1}$ and $B=-\cJ_2^{-1}\cJ_1$. Again, it is
straightforward to verify, using \eq{e:1}--\eq{e:7}, that
$\om^*=-\om$, $B^*=-B$, and that $\cJ_1=\om^{-1}B$,
$\cJ_2=-\om^{-1}$, $\cJ_3=\om+B\om^{-1}B$, $\cJ_4=-B\om^{-1}$.
Thus $E$ is $B$-symplectic, completing the proof.
\end{proof}

\subsection{Proofs for section
\ref{s:subquot}}\label{ss:proofssubquot} We begin with the
following lemma, whose proof is completely straightforward from
the definitions and is therefore omitted. Let $V$ be a real vector
space, let $E\subseteq V_{\bC}\oplus V^*_{\bC}$ be a maximally
isotropic subspace, and let $W\subseteq V$ be a subspace. Write
$\tau:V\oplus V^*\to V^*\oplus V$, $\tau_W:(V/W)\oplus (V/W)^*\to
(V/W)^*\oplus (V/W)$ for the isomorphisms interchanging the
summands. Define $E_W$ and $E_{V/W}$ by the formulas \eq{e:ew} and
\eq{e:evw}, respectively. Note that $\tau(E)$ is a maximally
isotropic subspace of $V^*_{\bC}\oplus V_{\bC}$ and $\Ann(W)$ is a
real subspace of $V^*$, so it makes sense to consider
$\tau(E)_{\Ann(W)}$. Finally, we identify $\Ann(W)$ with $V/W$ in
the obvious way.
\begin{lem}\label{l:duality}
With this identification, we have
\[
\tau(E)_{\Ann(W)}=\tau_W(E_{V/W}).
\]
\end{lem}
Our next result is a computation.
\begin{lem}\label{l:dim}
Let $V$ be a real vector space with a GC structure given by
$E\subseteq V_{\bC}\oplus V^*_{\bC}$. Let $W\subseteq V$ be a
subspace, and define $E_W$, $E_{V/W}$ by \eq{e:ew} and \eq{e:evw},
respectively. Then $\dim_{\bC} E_W=\dim_{\bR} W$ and $\dim_{\bC}
E_{V/W} = \dim_{\bR}(V/W)$.
\end{lem}
\begin{proof}
Note that the second equality is equivalent to the first, in view
of Lemma \ref{l:duality}. Thus we will only prove the first
equality. It is clear that $\rho(E_W)=\rho(E)\cap W_{\bC}$, so
projection to $W_{\bC}$ realizes $E_W$ as an extension
\[
0 \rar{} K \rar{} E_W \rar{} \rho(E)\cap W_{\bC} \rar{} 0,
\]
where $K$, in turn, fits into a short exact sequence
\[
0 \rar{} E\cap\Ann(W_{\bC}) \rar{} E\cap V_{\bC}^* \rar{} K \rar{}
0
\]
(the last nonzero arrow is given by restriction to $W$). Thanks to
the fact that $E$ is maximally isotropic, we have $E\cap
V_{\bC}^*=\Ann(\rho(E))$. Thus
\[
E\cap\Ann(W_{\bC})=E\cap
V_{\bC}^*\cap\Ann(W_{\bC})=\Ann(\rho(E)+W_{\bC}),
\]
which implies that
\[
\dim_{\bC}\bigl(E\cap\Ann(W_{\bC})\bigr)=\dim_{\bC} V_{\bC} -
\dim_{\bC}(\rho(E)+W_{\bC})=\dim_{\bR}V-\dim_{\bR}W -
\dim_{\bC}\rho(E)+\dim_{\bC}(\rho(E)\cap W_{\bC}).
\]
Therefore
\begin{eqnarray*}
\dim_{\bC}(E_W)&=&\dim_{\bC}(\rho(E)\cap W_{\bC})+\dim_{\bC}K \\
&=& \dim_{\bC}(\rho(E)\cap W_{\bC})+\dim_{\bC}(E\cap
V_{\bC}^*)-\dim_{\bC}\bigl(E\cap\Ann(W_{\bC})\bigr) \\
&=& \dimc(E\cap V_{\bC}^*)-\dimr V+\dimr W+\dimc\rho(E)=\dimr W,
\end{eqnarray*}
completing the proof.
\end{proof}

\begin{proof}[Proof of Proposition \ref{p:subspin}]
The product $f_1\wedge\dotsb\wedge f_k$ depends (up to a sign)
only on the span of the $f_i$'s, and we know from Theorem
\ref{t:spinors}(\ref{i:pure}) that the span must equal to $E\cap
V_{\bC}^*$. Thus, we are free to choose the $f_i$'s to be any
basis of this space. Since $E$ is maximally isotropic, we have
$E\cap V_{\bC}^*=\Ann(\rho(E))$, and now the existence of the
basis we are looking for follows from the short exact sequence
\[
0 \rar{} \Ann(\rho(E)+W_{\bC}) \rar{} \Ann(\rho(E)) \rar{j^*}
\Ann(\rho(E)\cap W_{\bC}) \rar{} 0,
\]
where we view the last space as a subspace in $W_{\bC}^*$. [The
surjectivity follows, for instance, from dimension counting.]

Now define $\phi_W$ by the formula given in Proposition
\ref{p:subspin}. Theorem \ref{t:spinors}(\ref{i:pure}) implies
that $\phi_W$ is a pure spinor, and \eq{e:ann} shows that
\[
E=E_{\phi}=\bigl\{ v-\iota_v(u)+f \,\big\vert\, v\in\rho(E),\ f\in
E\cap V_{\bC}^*=\Ann(\rho(E)) \bigr\},
\]
and the maximally isotropic subspace corresponding to $\phi_W$ is
\[
E_{\phi_W}=\bigl\{ w-\iota_w(u)+g \,\big\vert\, w\in\Ann(\Phi_W),\
g\in\Phi_W \bigr\},
\]
where $\Phi_W$ is the span of $j^*(f_1),\dotsc,j^*(f_l)$, i.e., by
construction, $\Phi_W=\Ann(\rho(E)\cap W_{\bC})\subseteq
W_{\bC}^*$. It is then clear that $E_W\subseteq E_{\phi_W}$, and
equality holds for reasons of dimension.
\end{proof}

\sbr

We now come to the proofs of Theorems \ref{t:classification} and
\ref{t:canon}. Let $V$ be an arbitrary GC vector space; we let
$\cJ\in\Aut_{\bR}(V\oplus V^*)$ denote the corresponding
orthogonal automorphism, and $E\subset V_{\bC}\oplus V^*_{\bC}$
the corresponding subspace. Let $C$ be the real subspace of $V$
such that $C_{\bC}=(E\cap V_{\bC})\oplus(\overline{E}\cap
V_{\bC})$. Then $C$ is obviously invariant under $\cJ$, so
$\cJ\restr{C}$ gives a canonical complex structure on $C$. It is
clear from Proposition \ref{p:graphgcsub}(a) that $C$, with this
structure, satisfies the graph condition. Furthermore, if $C$ is a
GC subspace of $V$, then Proposition \ref{p:graphgcsub}(c) implies
that the induced structure on $C$ is $\be$-complex, and the
underlying complex structure coincides with $\cJ\restr{C}$.

\sbr

On the other hand, there exists a real subspace $S\subseteq V$
such that $S_{\bC}=\rho(E)\cap\rho(\overline{E})$. It is easy to
check that $S$ is a GC subspace of $V$ (also see the proof of
Theorem \ref{t:classification} below). Moreover, Proposition
\ref{p:characterization}(c) implies that $S$ is the maximal GC
subspace of $V$ such that the induced structure is $B$-symplectic.
Next, $S_{\bC}=\Ann\bigl((E\cap V_{\bC}^*)\oplus (\overline{E}\cap
V_{\bC}^*)\bigr)$. If $\tilde{C}$ is the real subspace of $V^*$
such that $\tilde{C}_{\bC}=(E\cap V_{\bC}^*)\oplus
(\overline{E}\cap V_{\bC}^*)$, so that $\tilde{C}=\Ann(S)$, then,
applying the argument of the previous paragraph to the dual
structure $\tau(E)$, we see that $V/S\cong\tilde{C}^*$ acquires a
natural complex structure, and hence so does $\tilde{C}$.
Furthermore, if $V/S$ is a GC quotient of $V$, then by Proposition
\ref{p:subquot}, $\tilde{C}$ is a GC subspace of $V^*$, and now
the argument of the previous paragraph shows that $\tilde{C}$ is
$\be$-complex, whence $V/S$ is $B$-complex. On the other hand, if
we return to the situation of the previous paragraph and use our
results about $S$, we see that $V/C$ is always a GC quotient of
$V$, and the induced structure is $\be$-symplectic; and, moreover,
$C$ is the smallest subspace of $V$ with this property. This
completes the proof of Theorem \ref{t:canon}.

\sbr

To prove Theorem \ref{t:classification}, we begin by constructing
the maximal $B$-symplectic GC subspace $S\subseteq V$ as in the
previous paragraph. This characterization of $S$ shows that $S$
does not change if we make a $B$-field transform of the GCS on $V$
(and neither does the result of the theorem). Therefore, we
assume, from the beginning, that the pure spinor defining our GCS
has the form $\phi=\exp(i\Om)\wedge f_1\wedge\dotsb\wedge f_k$,
where $\Om\in\Wedge^2 V^*$ is a real $2$-form. We know that
$f_1,\dotsc,f_k$ form a basis for $E\cap V_{\bC}^*=\Ann(\rho(E))$,
and hence $f_1,\dotsc,f_k,\bar{f}_1,\dotsc,\bar{f}_k$ form a basis
for $\Ann(S_{\bC})$. By Theorem \ref{t:spinors}(\ref{i:pair}), we
have, up to a nonzero scalar multiple,
\[
0\neq \pair{\phi}{\bar{\phi}}_M = \Om^p\wedge
f_1\wedge\dotsb\wedge f_k\wedge
\bar{f}_1\wedge\dotsb\wedge\bar{f}_k,
\]
where $p=\dimr S$. This implies that the restriction of $\Om$ to
$S$ must be nondegenerate; write $\om=\Om\restr{S}$. Basic linear
algebra shows that if $W$ is the orthogonal complement to $S$ with
respect $\Om$, then $V=S\oplus W$ (direct sum of vector spaces).
Let $\pi_S:V\to S$, $\pi_W:V\to W$ be the two projections, and let
$g_j=f_j\restr{W_{\bC}}$ ($1\leq j\leq k$),
$u=i\Om\restr{W_{\bC}}$. Note that the $g_j,\bar{g}_j$ form a
basis of $W_{\bC}^*$. Now Proposition \ref{p:subspin} implies that
the pure spinors corresponding to $S$ and $W$ are given by
\[
\phi_S=\exp(i\om) \quad \text{and} \quad \phi_W=\exp(u)\wedge
f_1\wedge\dotsb\wedge f_k,
\]
respectively. Then Theorem \ref{t:spinors}(\ref{i:pair}) gives (up
to a nonzero scalar multiple)
\[
\pair{\phi_S}{\bar{\phi}_S}_M=\om^p \quad \text{and} \quad
\pair{\phi_W}{\bar{\phi}_W}_M = g_1\wedge\dotsb\wedge g_k\wedge
\bar{g}_1\wedge\dotsb\wedge\bar{g}_k;
\]
in particular, both $S$ and $W$ are GC subspaces of $V$. Moreover,
by construction, $W_{\bC}^*=(E_W\cap W_{\bC}^*)\oplus
(\overline{E}_W\cap W_{\bC}^*)$, whence $W$ is $B$-complex by
Proposition \ref{p:characterization}(a). Finally, since the $f_j$
annihilate $S_{\bC}$, and since $S$ is orthogonal to $W$ with
respect to $\Om$, it is easy to see that
\[
\phi=\pi_S^*(\phi_S)\wedge\pi_W^*(\phi_W),
\]
which shows that $V=S\oplus W$ as GC vector spaces, and completes
the proof of Theorem \ref{t:classification}.

\subsection{Proofs for section
\ref{s:split}}\label{ss:proofssplit} First let us briefly discuss
the proof of Proposition \ref{p:split}. It is obvious that
$N\oplus\Ann(N)$ is the orthogonal complement of $W\oplus\Ann(W)$
with respect to our standard pairing $\pairing$; since $\cJ$ is
orthogonal with respect to this pairing, we see that
$N\oplus\Ann(W)$ is invariant under $\cJ$. Now let
$\psi:W\oplus\Ann(N)\to W\oplus W^*$ and $\psi':N\oplus\Ann(W)\to
N\oplus N^*$ be the natural isomorphisms. Invariance under $\cJ$
implies that
\[
W_{\bC}\oplus\Ann(N_{\bC}) = \Bigl(E\cap
\bigl(W_{\bC}\oplus\Ann(N_{\bC})\bigr)\Bigr) \oplus
\Bigl(\overline{E}\cap\bigl(W_{\bC}\oplus\Ann(N_{\bC})\bigr)\Bigr);
\]
on the other hand, we have (by the definition of $E_W$)
\[
\psi\Bigl(E\cap \bigl(W_{\bC}\oplus\Ann(N_{\bC})\bigr)\Bigr)
\subseteq E_W,
\]
and similarly for $\overline{E}$. Comparing the dimensions, we see
that equality in the formula above holds. Since $\psi$ is an
isomorphism, it follows that $E_W\cap\overline{E}_W=(0)$. The rest
is obvious. We can now give the

\begin{proof}[Proofs of Propositions \ref{p:splitgraph}--\ref{p:graphgcsub}]
Consider a real vector space $V$ and a subspace $W\subseteq V$,
and suppose that $V$ and $W$ are equipped with GC structures $\cJ$
and $\cK$, respectively, written in the usual form as
\[
\cJ=\matr{\cJ_1}{\cJ_2}{\cJ_3}{\cJ_4} \text{ and }
\cK=\matr{\cK_1}{\cK_2}{\cK_3}{\cK_4}.
\]
The graph of the inclusion $W\hookrightarrow V$ is the subset
$\Ga\subset W\oplus V$ consisting of all pairs $(w,w)$ with $w\in
W$, while $\Ann(\Ga)$ is the subset of $W^*\oplus V^*$ consisting
of all pairs $(f,g)$ with $f=-g\restr{W}$. Let $\cJ_{tp}$ denote
the twisted product structure on $W\oplus V$; thus $\Ga$ is a
generalized isotropic subspace of $W\oplus V$ with respect to this
structure if and only if
$\cJ_{tp}(\Ga)\subseteq\Ga\oplus\Ann(\Ga)$. Explicitly, if $w\in
W$, so that $(w,w)\in\Ga$, then
\[
\cJ_{tp}(w,w,0,0)=\bigl(\cK_1(w),\cJ_1(w),-\cK_3(w),\cJ_3(w)\bigr).
\]
Therefore, $W$ (with the structure $\cK$) satisfies the graph
condition if and only if
\begin{equation}\label{e:nascgr}
\cK_1=\cJ_1\restr{W} \text{ (in particular, } W \text{ is
invariant under } \cJ_1 \text{), and } \cK_3(w)=\cJ_3(w)\restr{W}\
\forall\,w\in W.
\end{equation}
Note that \eq{e:nascgr} and \eq{e:5} imply that
$\cK_4(f\restr{W})=\cJ_4(f)\restr{W}$ for all $f\in V^*$. Thus,
the requirement that $W$ satisfies the graph condition with
respect to $\cK$ determines $\cK_1$, $\cK_3$ and $\cK_4$ uniquely.
In view of Lemma \ref{l:beta}, whose proof is given below, we see
that the rest of Propositions
\ref{p:splitgraph}--\ref{p:graphgcsub} follow from the following

\sbr

\noindent \textbf{Claim.} If $W$ is a GC subspace of $V$ which is
invariant under $\cJ_1$, and $\cK$ denotes the induced structure,
then the conditions \eq{e:nascgr} are satisfied.

To prove this claim, let $E$ be the $+i$-eigenspace of $\cJ$ in
$V_{\bC}$, and let $w\in W_{\bC}$. We can write $w=e+e'$, where
$e\in E$, $e'\in\overline{E}$. Then $\cJ(w)=i\cdot(e-e')$. Now
$\rho(e)+\rho(e')=w\in W_{\bC}$, and
$i\cdot(\rho(e)-\rho(e'))=\rho(\cJ(w))=\cJ_1(w)\in W_{\bC}$ by
assumption, whence $\rho(e),\rho(e')\in W_{\bC}$. Thus, in fact,
$(w,0)=(\rho(e),\rho^*(e)\restr{W_{\bC}})+(\rho(e'),\rho^*(e')\restr{W_{\bC}})$,
with $(\rho(e),\rho^*(e)\restr{W_{\bC}})\in E_W$ and
$(\rho(e'),\rho^*(e')\restr{W_{\bC}})\in\overline{E}_W$. Therefore
\[
\bigl(\cK_1(w),\cK_3(w)\bigr)=\cK(w,0)=i\cdot
(\rho(e),\rho^*(e)\restr{W_{\bC}}) - i\cdot
(\rho(e'),\rho^*(e')\restr{W_{\bC}})=\bigl(\cJ_1(w),\cJ_3(w)\restr{W_{\bC}}\bigr),
\]
completing the proof.
\end{proof}

\begin{proof}[Proof of Lemma \ref{l:beta}]
All the computations below that we omit are straightforward if one
uses \eq{e:1}--\eq{e:7}. First one checks that the condition that
$\be\in\Wedge^2 V$ transforms $\cJ$ into $\cJ'$ is equivalent to
the following equations:
\begin{equation}\label{e:be}
\be\cJ_3=\cJ_3\be=0, \quad -\cJ_1\be+\be\cJ_4=\cJ_2'-\cJ_2.
\end{equation}
Let $\cL=(1/2)\cdot(\cJ_2-\cJ_2')$. If $\cJ_1\be$ is to be
skew-symmetric, then (since $\be$ is also skew-symmetric) we must
have $\be\cJ_4=\be^*\cJ_1^*=(\cJ_1\be)^*=-\cJ_1\be$, which forces
$\cJ_1\be=\cL$. This equation and $\cJ_3\be=0$ imply that
$-\be=(\cJ_1^2+\cJ_2\cJ_3)\circ\be=\cJ_1\cL$, whence
$\be=(1/2)\cdot\cJ_1\circ(\cJ_2'-\cJ_2)$, proving the uniqueness
of $\be$. Conversely, if we define $\be$ by this formula, then one
easily checks that $\be^*=-\be$, and the equations \eq{e:be} are
satisfied. This completes the proof.
\end{proof}

\subsection{Proofs for section \ref{s:gcm}}\label{ss:proofsgcm}
\begin{proof}[Proof of Proposition \ref{p:cour}]
Part (a) is immediate from the equality $d(f\cdot g)+f\cdot
Y+g\cdot X=f\cdot (df+Y)+g\cdot (df+X)$. For part (b), let $(f,X)$
be a section of $\Cou_M$ defined near a point $m\in M$. Suppose
that $f(m)\neq 0$, then there exists an open set $U\ni m$ such
that $(f,X)$ is defined on $U$, and is nonvanishing there. Then it
is easy to check that $(f^{-1},-f^{-2}\cdot X)\in\Cou_M(U)$, and
$(f,X)\cdot (f^{-1},-f^{-2}\cdot X)=(1,0)$. Finally, part (c)
follows from the definition of integrability in terms of the
Courant bracket, and the easily verified identity
\[
\cou{df+X}{dg+Y}=\frac{1}{2}\bigl( d(X(g))-d(Y(f)) \bigr) + [X,Y]
\]
for all local sections $(f,X)$ and $(g,Y)$ of $C^\infty_M\oplus
T_{\bC}M$. We leave to the reader the straightforward but tedious
verification that $\{\cdot,\cdot\}$ satisfies the Leibnitz rule
and the definition of a Poisson algebra; we remark that here the
fact that $\overline{E}$ is isotropic is used in an essential way.
(In other words, $\{\cdot,\cdot\}$ doesn't have these properties
on the bigger sheaf $C^\infty_M\oplus T_{\bC}M$.)
\end{proof}

\begin{proof}[Proof of Theorem \ref{t:pullcour}]
The given map is clearly an ring homomorphism with respect to the
ring operations defined in Proposition \ref{p:cour}(a). Thus we
only need to check that this map takes $\Cou_M(U)$ into
$\Cou_S(U\cap S)$. Let $j:S\hookrightarrow M$ denote the inclusion
map. We'll show, more generally, that if $X+\xi$ is a section of
$\overline{E}$ over $U$, where $X\in\Ga(U,T_{\bC}M)$ and
$\xi\in\Ga(U,T^*_{\bC}M)$, then $\pi(X\restr{S})+j^*(\xi)$ is a
section of $\overline{E}_S$ over $U\cap S$. To that end, write, on
$U\cap S$: $X=Y+Z$, $\xi=\eta+\zeta$, where $Y$, $Z$, $\eta$ and
$\zeta$ are sections over $U\cap S$ of $T_{\bC}S$, $N_{\bC}$,
$\Ann(N_{\bC})\subset T_{\bC}^*M\restr{S}$ and
$\Ann(T_{\bC}S)\subset T_{\bC}^*M\restr{S}$, respectively. It is
clear that $Y=\pi(X\restr{S})$ and
$\eta\restr{T_{\bC}S}=j^*(\xi)$. Thus it remains to observe that
$(Y,\eta)$ is a section of $\overline{E}$. This follows from the
decomposition
\[
X+\xi=(Y+\eta)+(Z+\zeta)\in \bigl(T_{\bC}S+\Ann(N_{\bC})\bigr)
\oplus \bigl(N_{\bC}+\Ann(T_{\bC}S)\bigr)
\]
and the fact that each of the two direct summands on the right is
$\cJ$-invariant, by the definition of a split submanifold.
\end{proof}

\subsection{Proofs for section
\ref{s:category}}\label{ss:proofscategory}
\begin{proof}[Proof of Theorem \ref{t:graphiso}]
The theorem follows immediately from its linear counterpart
(Theorem \ref{it:4} in the introduction), which we find convenient
to restate as follows. Let $V$ be a real vector space, and let
$\cJ$, $\cK$ be two GC structures on $V$. If the diagonal
$\De\subset V\oplus V$ is generalized Lagrangian with respect to
the direct sum of the structures $\cJt$ and $\cK$, then $\cJ=\cK$.

\sbr

Now $\Ann(\De)=\{(f,-f)\in V^*\oplus V^*\}$, and if $v\in V$,
$f\in V^*$, then
\[
(\cJt\oplus\cK)(v,v,f,-f)= \bigl(
\cJ_1(v)-\cJ_2(f),\cK_1(v)-\cK_2(f),-\cJ_3(v)+\cJ_4(f),\cK_3(v)-\cK_4(f)
\bigr).
\]
We see that $\De\oplus\Ann(\De)$ is invariant under
$\cJt\oplus\cK$ if and only if $\cJ_1=\cK_1$, $\cJ_2=\cK_2$,
$\cJ_3=\cK_3$ and $\cJ_4$=$\cK_4$, completing the proof.
\end{proof}

\begin{proof}[Proof of Theorem \ref{t:closcomp}]
In the following proof, we will always implicitly identify
$(V\oplus W)^*$ with $V^*\oplus W^*$, and similarly for $W\oplus
Z$, $V\oplus Z$. We begin with the following

\sbr

\noindent \textbf{Claim.} Let $(f,h)\in V^*\oplus Z^*$. Then
$(f,h)\in\Ann(\Phi\circ\Ga)$ if and only if there exists $g\in
W^*$ such that $(f,g)\in\Ann(\Ga)$ and $(-g,h)\in\Ann(\Phi)$.

Indeed, suppose first that there is a $g\in W^*$ with
$(f,g)\in\Ann(\Ga)$, $(-g,h)\in\Ann(\Phi)$. If
$(v,z)\in\Phi\circ\Ga$, then there exists $w\in W$ such that
$(v,w)\in\Ga$ and $(w,z)\in\Phi$. Then $f(v)+g(w)=0$ and
$-g(w)+h(z)=0$, whence $f(v)+h(z)=0$. This shows that
$(f,h)\in\Ann(\Phi\circ\Ga)$.

Conversely, let $(f,h)\in\Ann(\Phi\circ\Ga)$. Let $\La$ (resp.,
$\Psi$) denote the projection of $\Ga$ (resp., $\Phi$) onto $W$.
Note that if $(v,0)\in\Ga$, then $(v,0)\in\Phi\circ\Ga$, so that
$f(v)=-h(0)=0$; similarly, if $(0,z)\in\Phi$, then
$(0,z)\in\Phi\circ\Ga$, whence $h(z)=0$. This shows that we can
define a linear functional $g_1$ on $\La$ by the rule
$g_1(w)=-f(v,w)$ for any $v\in V$ such that $(v,w)\in\Ga$.
Similarly, we can define a linear functional $g_2$ on $\Psi$ by
the rule $g_2(w)=h(w,z)$ for any $z\in Z$ such that
$(w,z)\in\Phi$. Then $g_1$ and $g_2$ agree on $\La\cap\Psi$,
because if $w\in\La\cap\Psi$, then $(v,w)\in\Ga$ and
$(w,z)\in\Phi$ for some $v\in V$, $z\in Z$, whence
$g_1(w)=-f(v)=h(z)=g_2(w)$. Thus, there exists a linear functional
$g\in W^*$ such that $g\restr{\La}=g_1$ and $g\restr{\Psi}=g_2$.
By construction, we have $(f,g)\in\Ann(\Ga)$ and
$(g,h)\in\Ann(\Psi)$, proving the claim.

\sbr

Now let us prove part (a) of the theorem. Let $\cJ_V$, $\cJ_W$ and
$\cJ_Z$ denote the automorphisms defining the GC structures on
$V$, $W$ and $Z$, respectively. We will write them in the usual
$2\times 2$ matrix form. Suppose that $\Ga\subset V\oplus W$ and
$\Phi\subset W\oplus Z$ are generalized isotropic subspace with
respect to the twisted product structures. If
$(v,z)\in\Phi\circ\Ga$, then there exists $w\in W$ with
$(v,w)\in\Ga$ and $(w,z)\in\Phi$. Then we get
\[
(\cJt_V\oplus\cJ_Z)(v,z,0,0)=\bigl( \cJ_{V1}(v), \cJ_{Z1}(z),
-\cJ_{V3}(v), \cJ_{Z3}(z) \bigr).
\]
We have
\[
(\cJt_V\oplus\cJ_W)(v,w,0,0)=\bigl( \cJ_{V1}(v), \cJ_{W1}(w),
-\cJ_{V3}(v), \cJ_{W3}(w) \bigr)
\]
and
\[
(\cJt_W\oplus\cJ_Z)(w,z,0,0)=\bigl( \cJ_{W1}(w), \cJ_{Z1}(z),
-\cJ_{W3}(w), \cJ_{Z3}(z) \bigr),
\]
whence $\bigl(\cJ_{V1}(v),\cJ_{W1}(w)\bigr)\in\Ga$,
$\bigl(\cJ_{W1}(w),\cJ_{Z1}(z)\bigr)\in\Phi$,
$\bigl(-\cJ_{V3}(v),\cJ_{W3}(w)\bigr)\in\Ann(\Ga)$, and
$\bigl(-\cJ_{W3}(w),\cJ_{Z3}(z)\bigr)\in\Ann(\Phi)$, by
assumption. This implies that
$\bigl(\cJ_{V1}(v),\cJ_{Z1}(z)\bigr)\in\Phi\circ\Ga$, and,
according to the claim above,
$\bigl(-\cJ_{V3}(v),\cJ_{Z3}(z)\bigr)\in\Ann(\Phi\circ\Ga)$. This
proves that $\Phi\circ\Ga$ is generalized isotropic with respect
to the twisted product structure on $V\oplus Z$.

\sbr

For part (b), we assume that $\Ga$ and $\Phi$ are generalized
coisotropic with respect to the twisted product structures. Then
we start with $(f,h)\in\Ann(\Phi\circ\Ga)$; according to the claim
above, we can find a $g\in W^*$ such that $(f,g)\in\Ann(\Ga)$ and
$(-g,h)\in\Ann(\Phi)$. Then we proceed exactly as in the previous
paragraph, and we deduce that
\[
(\cJt_V\oplus\cJ_Z)(0,0,f,h)\in
(\Phi\circ\Ga)\oplus\Ann(\Phi\circ\Ga).
\]

\sbr

Finally, part (c) follows from parts (a) and (b).
\end{proof}

\end{document}